\documentstyle[11pt,leqno]{article}

\include{psfig}
\include{epsf}

\newtheorem{theorem}{Theorem}[section]

\newtheorem{proposition}[theorem]{Proposition}

\newtheorem{remark}{Remark}[section]
\newtheorem{example}{Example}[section]

\catcode`\@=11
\renewcommand{\section}{
        \setcounter{equation}{0}
        \@startsection {section}{1}{\z@}{-3.5ex plus -1ex minus
        -.2ex}{2.3ex plus .2ex}{\Large\bf}
}
\catcode`\@=12

\def\reals{{\rm\vrule depth0ex width.4pt\kern-.08em R}}
\def\bbbz{{\mathchoice {\hbox{$\sf\textstyle Z\kern-0.4em Z$}}
{\hbox{$\sf\textstyle Z\kern-0.4em Z$}}
{\hbox{$\sf\scriptstyle Z\kern-0.3em Z$}}
{\hbox{$\sf\scriptscriptstyle Z\kern-0.2em Z$}}}}

\newcommand{\nc}{\newcommand}
\nc{\W}{{\bf W}}
\nc{\A}{{\bf A}}
\nc{\bL}{{\bf L}}
\nc{\bH}{{\bf H}}
\nc{\C}{{\cal C}}

\def\eq#1{(\ref{e:#1})}

\def\elabel#1{\label{e:#1}}

\begin{document}
\begin{center}
\Large \bf Product-form solutions for integrated services packet
networks and cloud computing systems~\footnote{Partial results,
graphes, and tables are briefly summarized and reported in MICAI'06.
This enhanced version with the latest discussions and complete
proofs of results is a journal version of the short conference
report.}
\end{center}
\begin{center}
Wanyang Dai~\footnote{Supported by National Natural Science
Foundation of China with Grant No.10371053, Grant No.10971249, Grant
No.11371010.}\\
Department of Mathematics and State Key Laboratory of Novel
Software Technology\\
Nanjing University\\
Nanjing 210093\\
China\\
Email: nan5lu8@netra.nju.edu.cn\\
Submitted: 3 September 2013
\end{center}

\vskip 0.1 in
\begin{abstract}

We iteratively derive the product-form solutions of stationary
distributions of priority multiclass queueing networks with
multi-sever stations. The networks are Markovian with exponential
interarrival and service time distributions. These solutions can be
used to conduct performance analysis or as comparison criteria for
approximation and simulation studies of large scale networks with
multi-processor shared-memory switches and cloud computing systems
with parallel-server stations. Numerical comparisons with existing
Brownian approximating model are provided
to indicate the effectiveness of our algorithm.\\


\end{abstract}

\section{Introduction}\label{introduction}

\begin{figure}[tbh]
\centerline{\epsfxsize=4.5in\epsfbox{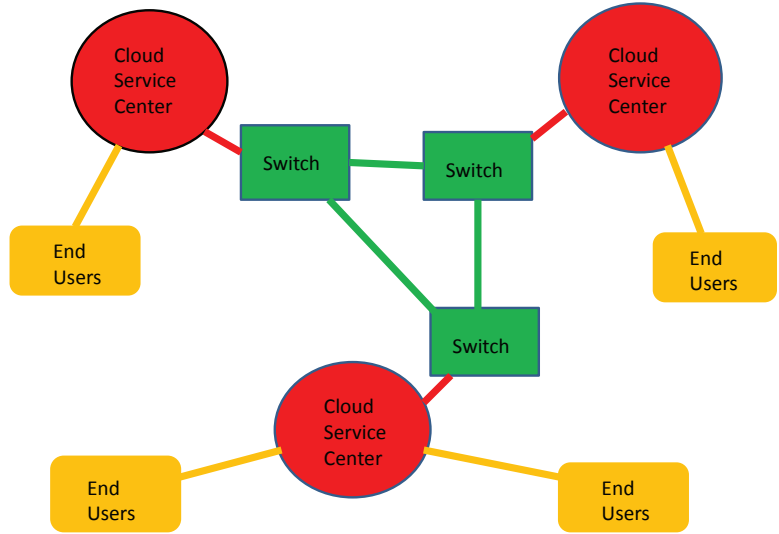}} \caption{An
integrated services cloud computing based network} \label{cloudnet}
\end{figure}

At present, integrated services packet networks (ISPN) are widely
used to transport a wide range of information such as voice, video
and data. It is foreseeable that this integrated services pattern
will be one of the major techniques in the future cloud computing
based communication systems and Internet. The introduction of
concept and architecture of cloud computing can be found in the
white paper of Sun-microsystems~\cite{sun:intclo}. A possible cloud
computing based telecommunication network architecture (i.e., a
large-scale network infrastructure as a service) is designed in
Figure~\ref{cloudnet}, where an end-user may require service (or
services) from single local cloud computing center or multiple local
and remote cloud computing centers. Among these centers, they
communicate each other by using core switching network system (note
that the switching system itself can also be independently viewed
and handled as a cloud computing system with multiple service
pools).

The speed and efficiency of core switching network systems are the
bottleneck in realizing high-speed owing to the drastic improvement
in transmission speed and reliability of optical fiber. In an ISPN
network, information is partitioned into equal or non-equal size
packets depending on the employed protocol such as Internet protocol
(IP), asynchronous transfer mode (ATM), and frame relay (FR). For
the purpose of transmission, each packet consists of user's data
payload to be transmitted, a header containing control information
(e.g., source and destination addresses, packet type, priority,
etc.), and a checksum used for error control. The high-speed ISPN
networks require fast packet switches to move packets along their
respective virtual paths. The switches are computers with processing
and storage capabilities. The main objective of a switch is to route
an incoming packet arriving along a particular input link to a
specified output link. More precisely, once the incoming packet is
entirely received, the switch examines its header and selects the
next outgoing link. In other words, packets are transmitted in a
{\em store-and-forward} manner.

Various information can be classified into a fixed number of
different types with separate requirements of quality of service
(e.g., different end-to-end delays and packet loss ratios). Real
time traffic packets with stringent delay constraint (e.g.,
interactive audio/video) are endowed service priority. In the
meanwhile, it is imperative to size buffers for non-real time
traffic packets (e.g., data), which can tolerate higher delay but
demand much smaller packet loss ratios. Hence, efficient switching
and buffer management schemes are needed for switches.

Currently, three basic techniques (see, e.g., Tobagi
\cite{tob:faspac}) are designed to realize the switching function:
{\em space-division, shared-medium, and shared-memory}. In our
switches, the shared-memory technique is employed, which are
comprise of a single dual-ported memory shared by all input and
output lines. Packets arriving along input lines are multiplexed
into a single stream that is fed to the common memory for storage.
Inside the memory, packets are organized into different output
buffers, one for each output line. In each output buffer, memory can
be further divided into priority output queues, one for each type of
traffic. In the meantime, an output stream of packets is formed by
retrieving packets from the output queues sequentially, one per
queue. Among different traffic types, packets are retrieved
according to their priorities. Inside each type, packets are
retrieved under the first-in first-out (FIFO) service discipline.
Then, the output stream is demultiplexed and packets are transmitted
on the output lines. There are some other ways to deal with the
output queues such as processor-sharing among output lines (see,
e.g., Tobagi \cite{tob:faspac}, Arpaci and Copeland
\cite{arpcop:bufman}). We will address these issues elsewhere (see,
e.g., Dai~\cite{dai:optrat}).

In addition to switching, queueing is another main functionality of
packet switches. The introduction of queueing function is owing to
multiple packets arriving simultaneously from separate input lines
and owing to the randomness of packet arrivals from both outside and
internal routing of the network. There are three possibilities for
queueing in a switch: {\em input-queueing, output-queueing}, and
{\em shared-memory-queueing} (see, e.g., Schwarz \cite{sch:broint}).
Shared-memory-queueing mechanism is employed in our switches since
it has numerous advantages with respect to switching and queueing
over other schemes. For example, both switching and queueing
functions can be implemented together by controlling the memory read
and write properly. Furthermore, modifying the memory read and write
control circuit makes the shared-memory switch sufficiently flexible
to perform functions such as priority control and other operations
(see, e.g., Endo {\em et al.}~\cite{endkoz:shabuf}). Thus,
shared-memory switches gained popularity among switch vendors.

\begin{figure}
\begin{picture}(100,200)(10,20)
\put(131,180){switch 1} \put(150,60){\circle{20}{$c_{1}$}}
\put(150,120){\circle{20}{2}} \put(150,145){\circle{20}{1}}
\put(150,80){\circle*{5}} \put(150,90){\circle*{5}}
\put(150,100){\circle*{5}} \put(130,35){\line(0,1){135}}
\put(170,35){\line(0,1){135}} \put(130,35){\line(1,0){40}}
\put(130,170){\line(1,0){40}} \put(160,60){\vector(3,2){37}}
\put(158,115){\vector(2,1){38}} \put(158,60){\vector(2,-1){50}}
\put(210,32){exit} \put(158,150){\vector(2,-1){32}}
\put(194,235){\vector(0,-1){98}{external arrivals}}
\put(60,80){\line(1,0){50}} \put(60,90){\line(1,0){50}}
\put(110,80){\line(0,1){10}} \put(105,80){\line(0,1){10}}
\put(100,80){\line(0,1){10}} \put(95,80){\line(0,1){10}}
\put(90,80){\line(0,1){10}} \put(85,80){\line(0,1){10}}
\put(80,80){\line(0,1){10}} \put(75,80){\line(0,1){10}}
\put(20,85){\vector(1,0){40}} \put(70,95){class 3}
\put(110,85){\vector(4,-3){30}} \put(410,126){exit}
\put(394,190){feedback}

\put(60,130){\line(1,0){50}} \put(60,140){\line(1,0){50}}
\put(110,130){\line(0,1){10}} \put(105,130){\line(0,1){10}}
\put(100,130){\line(0,1){10}} \put(20,135){\vector(1,0){40}}
\put(70,145){class 1} \put(110,135){\vector(3,-1){30}}
\put(110,135){\vector(3,1){30}}
\put(270,180){switch 2} \put(290,60){\circle{20}{$c_{2}$}}
\put(290,120){\circle{20}{2}} \put(290,145){\circle{20}{1}}
\put(290,80){\circle*{5}} \put(290,90){\circle*{5}}
\put(290,100){\circle*{5}} \put(270,35){\line(0,1){135}}
\put(310,35){\line(0,1){135}} \put(270,35){\line(1,0){40}}
\put(270,170){\line(1,0){40}} \put(298,115){\vector(2,1){38}}
\put(298,150){\vector(2,-1){32}} \put(290,100){\vector(2,1){70}}
\put(300,60){\vector(3,2){39}}

\put(200,80){\line(1,0){50}} \put(200,90){\line(1,0){50}}
\put(250,80){\line(0,1){10}} \put(245,80){\line(0,1){10}}
\put(240,80){\line(0,1){10}} \put(235,80){\line(0,1){10}}
\put(230,80){\line(0,1){10}} \put(225,80){\line(0,1){10}}
\put(220,80){\line(0,1){10}} \put(215,80){\line(0,1){10}}
\put(330,85){\vector(1,0){100}} \put(410,87){exit}
\put(210,95){class 4} \put(390,85){\vector(0,-1){65}}
\put(390,20){\vector(-1,0){350}} \put(394,22){feedback}
\put(40,20){\vector(0,1){65}} \put(250,85){\vector(4,-3){30}}

\put(200,130){\line(1,0){50}} \put(200,140){\line(1,0){50}}
\put(250,130){\line(0,1){10}} \put(245,130){\line(0,1){10}}
\put(240,130){\line(0,1){10}} \put(235,130){\line(0,1){10}}
\put(230,130){\line(0,1){10}} \put(330,135){\vector(1,0){100}}
\put(210,145){class 2} \put(390,135){\vector(0,1){65}}
\put(390,200){\vector(-1,0){350}} \put(40,200){\vector(0,-1){65}}
\put(250,135){\vector(3,-1){30}} \put(250,135){\vector(3,1){30}}
\put(250,135){\vector(1,-1){35}}
\end{picture}
\caption{\small A two-type and four-class network with two
multiprocessor-switches. Type 1 includes class 1 and class 2. Type 2
consists of class 3 and class 4. Type 1 owns the service priority.}
\label{fig:1}
\end{figure}
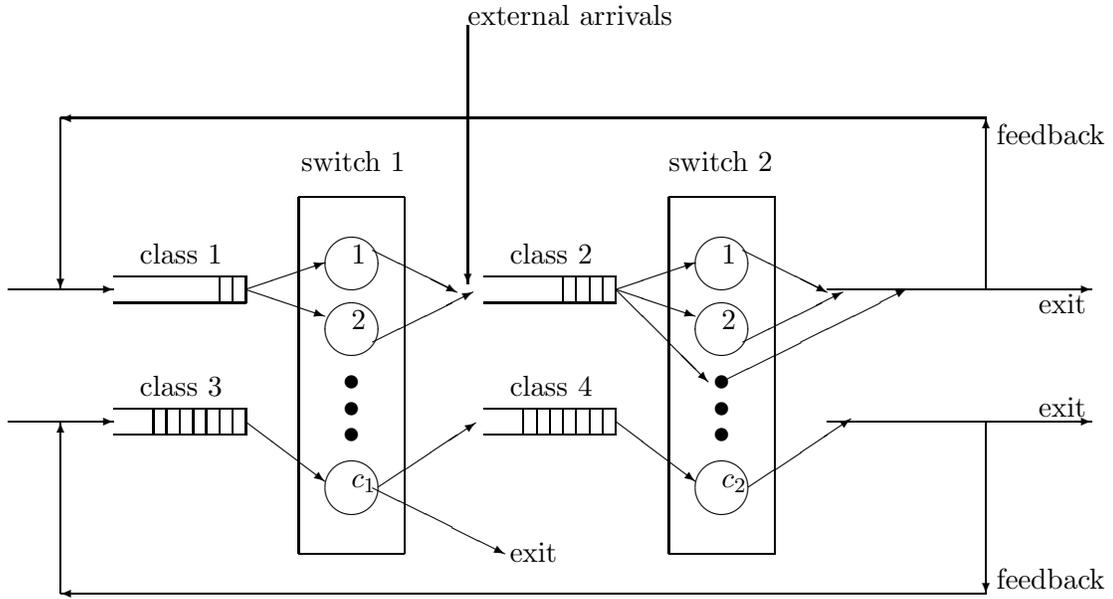

\begin{figure*}
\begin{picture}(100,200)(10,20)
\put(120,190){\line(1,0){230}} \put(120,50){\line(1,0){230}}
\put(120,50){\line(0,1){140}} \put(350,50){\line(0,1){140}}

\put(175,170){\line(1,0){40}} \put(175,130){\line(1,0){40}}
\put(175,170){\line(0,-1){40}} \put(215,170){\line(0,-1){40}}
\put(70,160){\line(1,0){105}} \put(70,140){\line(1,0){105}}
\put(215,160){\line(1,0){40}} \put(215,140){\line(1,-1){40}}
\put(185,148){$2\times 2$}

\put(175,110){\line(1,0){40}} \put(175,70){\line(1,0){40}}
\put(175,110){\line(0,-1){40}} \put(215,110){\line(0,-1){40}}
\put(70,100){\line(1,0){105}} \put(70,80){\line(1,0){105}}
\put(215,100){\line(1,1){40}} \put(215,80){\line(1,0){40}}
\put(185,88){$2\times 2$}

\put(255,170){\line(1,0){40}} \put(255,130){\line(1,0){40}}
\put(255,170){\line(0,-1){40}} \put(295,170){\line(0,-1){40}}
\put(295,160){\line(1,0){105}} \put(295,140){\line(1,0){105}}
\put(265,148){$2\times 2$}

\put(255,110){\line(1,0){40}} \put(255,70){\line(1,0){40}}
\put(255,110){\line(0,-1){40}} \put(295,110){\line(0,-1){40}}
\put(295,100){\line(1,0){105}} \put(295,80){\line(1,0){105}}
\put(265,88){$2\times 2$}

\end{picture}
\caption{A two-stage switch with four dual-ported shared-memory
switching modules.} \label{fig:mstage}
\end{figure*}
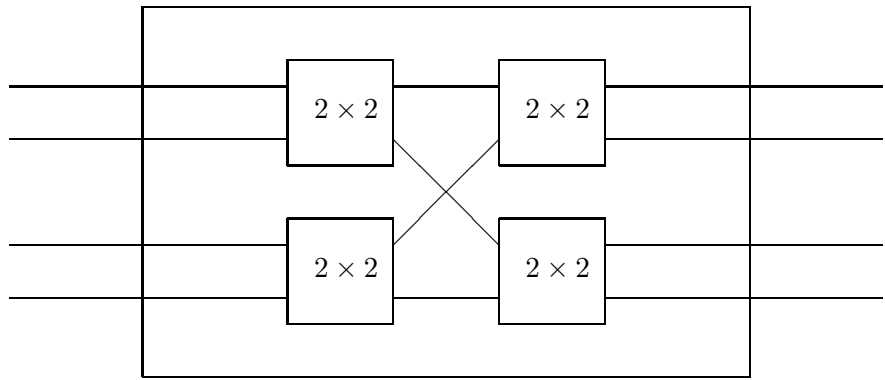

Nevertheless, no matter what technology is employed in implementing
the switch, it places certain limitations on the size of the switch
and line speeds. Currently, two possible ways can be used to build
large switches to match the transmission speed of optical fiber. The
first one is to adopt parallel processors (see, e.g.,
Figure~\ref{fig:1}).
The second one is to interconnect many switches (known as switch
modules) in a {\em multistage} configuration to build large switches
(see, e.g., Figure~\ref{fig:mstage}).
The remaining issue is how to reasonably allocate resources of these
switches and efficiently evaluate the system performance.

The statistical characteristics of packet interarrival times and
packet sizes have a major impact on switch hardware and software
designs owing to the consideration of network performance. How to
more effectively identify packet traffic patterns is a very active
and involved research field (see, e.g., Nikolaidis and Akyildiz
\cite{nikaky:ovesou} for a survey). Independent and identical
distribution (i.i.d.) is the popular assumption for these times and
packet sizes. Doubly stochastic renewal process introduced in
Dai~\cite{dai:optrat} is the latest definition and generalization
related to input traffic and service processes for a wireless
network under random environment. The effectiveness of these
characteristics is supported by recent discoveries in Cao {\em et
al.}~\cite{caocle:inttra} and Dai~\cite{dai:contru,dai:heatra}.

Note that, in all circumstances, it is imperative to find
product-form solutions for those queueing networks under suitable
conditions to conduct performance analysis or provide comparison
criteria to show the effectiveness of approximation and/or
simulation studies (see, e.g., Dai {\em et
al.}~\cite{dai:broapp,dai:difapp,daidai:heatra,shechedaidai:finele}).
Therefore, the aim of the current paper is to derive product-form
solutions iteratively for priority multiclass queueing networks with
multi-server stations when the interarrival and service times are
exponentially distributed. Our findings are new and our discussions
are different from those as summarized in Serfozo~\cite{ser:intsto}.
Furthermore, numerical comparisons with existing Brownian
approximating model are provided to indicate the effectiveness of
our algorithm.

The remainder of this paper is organized as follows. The open
priority multiclass queueing network associated with high-speed ISPN
is described in Section \ref{c:model}. Our main results including
product-form solutions and performance comparisons are presented in
Section \ref{distrib}. Numerical comparisons are given in
Section~\ref{comparison}. The proofs of our main theorems are
provided in Section~\ref{proof}. The conclusion of the paper is
presented in Section~\ref{conclusion}.

\section{The queueing network model}\label{c:model}

We consider a queueing network that has $J$ multiserver-stations.
Each station indexed by $j\in\{1,...,J\}$ owns $c_{j}$ servers and
has an infinite capacity waiting buffer. In the network, there are
$I$ job types. Each type consists of $J$ job classes that are
distributed at different stations. Therefore, the network is
populated by $K\;(=IJ)$ job classes which are labeled by
$k\in\{1,...,K\}$. Upon the arrival of a job of a type from outside
the network, it may only receive service for part of $J$ classes and
may visit a particular class more than once (but at most finite many
times). Then, it leaves the network (i.e., the network is open). At
any given time during its lifetime in the network, the job belongs
to one of the classes and changes classes at each time a service is
completed. All jobs within a class are served at a unique station
and more than one class might be served at a station (so-called {\em
multiclass} queueing network). The ordered sequence of classes that
a job visits in the network is named {\em a route}. Inter-routing
among different job types  is not allowed throughout the entire
network.

We use ${\cal C}(j)$ to denote the set of classes belonging to
station $j$. Let $s(k)$ denote the station to which class $k$
belongs. We implicitly set $j=s(k)$ when $j$ and $k$ appear
together. Associated with each class $k$, there are two i.i.d
sequences of random variables (r.v.), $u_{k}=\{u_{k}(i),i\geq 1\}$
and $v_{k}=\{v_{k}(i),i\geq 1\}$, an i.i.d sequence of
$K$-dimensional random vectors, $\phi^{k}=\{\phi^{k}(i),i\geq 1\}$,
and two real numbers, $\alpha_{k}\geq 0$ and $m_{k}=1/\mu_{k}>0$. We
suppose that the $3K$ sequences
\begin{eqnarray}
u_{1},...,u_{K},v_{1},...,v_{K},\phi^{1},...,\phi^{K}
\elabel{primitiveI}
\end{eqnarray}
are mutually independent. The initial r.v.s $u_{k}(1)$ and
$v_{k}(1)$ have means $\Huge{E}[u_{k}(1)]=1/\alpha_{k}$ and
$\Huge{E}[v_{k}(1)]=m_{k}$ respectively. For each $i\in\{1,2,...\}$,
$u_{k}(i)$ denotes the interarrival time between the $(i-1)$th and
the $i$th {\em externally} arrival job at class $k$. Furthermore,
$v_{k}(i)$ denotes the service time for the $i$th class $k$ job. In
addition, $\phi^{k}(i)$ denotes the routing vector of the $i$th
class $k$ job. We allow $\alpha_{k}=0$ for some classes $k\in{\cal
E}\equiv\{k:\alpha_{k}\neq 0\}$. Then, it follows that $\alpha_{k}$
and $\mu_{k}$ are the external arrival rate and service rate for
class $k$ respectively. We assume that the routing vector
$\phi^{k}(i)$ takes values in $\{e_{0},e_{1},...,e_{K}\}$, where
$e_{0}$ is the $K$-dimensional vector of all 0's. For $l=1,...,K$,
$e_{l}$ is the $K$-dimensional vector with $l$th component 1 and
other components 0. When $\phi^{k}(i)=e_{l}$, the $i$th job
departing class $k$ becomes a class $l$ job. Let
$p_{kl}=P\{\phi^{k}(i)=e_{l}\}$ be the probability that a job
departing class $k$ becomes a class $l$ job (of the same type).
Thus, the corresponding $K\times K$ matrix $P=(p_{kl})$ is {\em
routing matrix} of the network. Furthermore, the matrix
\begin{eqnarray}
{\cal Q}=I+P'+(P')^{2}+\cdot\cdot\cdot
\elabel{matrixQ}
\end{eqnarray}
is finite, i.e., $(I-P')$ is invertible with ${\cal Q}=(I-P')^{-1}$
since the network is open. The symbol $'$ on a vector or a matrix
denotes the transpose and $I$ denotes the identity matrix.

We use $\lambda_{k}$ for $k\in\{1,...,K\}$ to denote the overall
arrival rate to class $k$, including both external arrivals and
internal transitions. Then, we have the following {\em traffic
equation}
\begin{eqnarray}
\lambda_{k}=\alpha_{k}+\sum_{l=1}^{K}\lambda_{l}p_{lk},
\elabel{traequ}
\end{eqnarray}
or in its vector form, $\lambda=\alpha+P'\lambda$ (all vectors in
this paper are to be interpreted as column vectors unless explicitly
stated otherwise). Note that the unique solution $\lambda$ of
\eq{traequ} is given by $\lambda={\cal Q}\alpha$. For each
$\lambda_{k}$, if there is a long-run average rate of flow into the
class which equals to the long-run average rate out of that class,
this rate will equal $\lambda_{k}$. Furthermore, we define the {\em
traffic intensity} $\rho_{j}$ for station $j$ as follows
\begin{eqnarray}
\rho_{j}=\sum_{k\in{\cal C}(j)}
\frac{\lambda_{k}}{c_{j}\mu_{k}}.
\elabel{rhoj}
\end{eqnarray}
$\rho_{j}$ with $\rho_{j}\leq 1$ is also referred to as the nominal
fraction of time that station $j$ is non-idle.

The order of jobs being served at each station is dictated by a
service discipline. In the current research, we restrict our
attention to static buffer priority (SBP) service disciplines under
which the classes at each station are assigned a fixed rank (with no
ties). In our queueing network, each type of jobs is assigned the
same priority rank at every station where it possibly receives
service. When a server within a station switches from one job to
another, the new job will be taken from the leading (or longest
waiting) job at the highest rank non-empty class at the server's
station. Within each class, it is assumed that jobs are served on
the first-in first-out (FIFO) basis. We suppose that the discipline
employed is {\em non-idling}, i.e., a server is never idle when
there are jobs waiting to be served at its station. We also assume
that the discipline is {\em preemptive-resume}, i.e., when a job
arrives at a station with all servers busy and if the job is with a
higher rank than at least one of the jobs currently being served,
one of the lower rank job services is interrupted; when there is a
server available to the interrupted service, it continues from where
it left off. For convenience and without loss of generality, we use
consecutive numbers to index the classes that have the same priority
rank at stations $1$ to $J$. In other words, the highest priority
classes for station 1 to $J$ are indexed by $1$ to $J$, the second
highest priority classes are indexed by $J+1$ to $2J$ ,..., and the
lowest priority classes are indexed by $K-J+1$ to $K$. An example of
such a two-station network is given in Figure~\ref{fig:1}. In the
network, type 1 traffic possibly requires class 1 and class 2
services, type 2 traffic possibly requires class 3 and class 4
services, classes in type 1 have the higher priority at their
corresponding stations.

Finally, we define the cumulative arrival, cumulative service and
cumulative routing processes by the sums,
\begin{eqnarray}
U_{k}(n)=\sum_{i=1}^{n}u_{k}(i),\;\;\;
V_{k}(n)=\sum_{i=1}^{n}v_{k}(i),\;\;\;
\Phi^{k}(n)=\sum_{i=1}^{n}\phi^{k}(i),
\elabel{uvphi}
\end{eqnarray}
where the $l$th component $\Phi^{k}_{l}(n)$ of $\Phi^{k}(n)$
is the cumulative number of jobs to class $l$ for the
first $n$th jobs leaving class $k$ with $n=1,2,...$ and
$k,l\in\{1,...,K\}$. Then, we define
\begin{eqnarray}
&&E_{k}(t)\equiv\mbox{max}\{n\geq 0,U_{k}(n)\leq t\},
\elabel{cumE}\\
&&S_{k}(t)\equiv\mbox{max}\{n\geq 0,V_{k}(n)\leq t\},
\elabel{cumS}\\
&&A_{k}(t)\equiv E_{k}(t)+
\sum_{l=1}^{K}\Phi^{l}_{k}(S_{l}(t)).
\elabel{cumA}
\end{eqnarray}
Note that $E_{k}(t)$ denotes the total number of external arrivals
to class $k$ in time interval $[0,t]$. $S_{k}(t)$ represents the
total number of class $k$ jobs which have finished service in
$[0,t]$. $A_{k}(t)$ is the total arrivals to class $k$ in $[0,t]$
including both external arrivals and internal transitions.

\section{Steady-state queue length distributions}\label{distrib}

We use $Q_{k}(t)$ to denote the number of class $k$ jobs in station
$j=s(k)$ with $j\in\{1,2,...,J\}$ and $k\in\{1,2,...,K\}$ at time
$t$. It is called the {\em queue length process} for class $k$ jobs.
For convenience, let $Q_{(i,j)}(t)$ and $n_{(i,j)}$ be the
$(j-i+1)$-dimensional queue length process and $(j-i+1)$-dimensional
state vectors respectively. They are given by
\begin{eqnarray}
Q_{(i,j)}(t)=\left(Q_{i}(t),...,Q_{j}(t)\right),\;\;\;\;
n_{(i,j)}=\left(n_{i},...,n_{j}\right)
\elabel{qnij}
\end{eqnarray}
for $j>i$ and $i,j\in\{1,2,...,K\}$ and nonnegative integers
$n_{i},...,n_{j}$. Then, we use $P_{n_{(1,K)}}(t)$ to denote the
probability of system at state $n_{(1,K)}=(n_{1},...,n_{K})$ and at
time $t$, i.e.,
\begin{eqnarray}
P_{n_{(1,K)}}(t)=P\left\{Q_{(1,K)}(t)=n_{(1,K)}\right\}.
\elabel{P(t)}
\end{eqnarray}
Furthermore, let
\begin{eqnarray}
P_{n_{(1,K)}}=P\left\{Q_{(1,K)}(\infty)=n_{(1,K)}\right\}
\elabel{steadyP}
\end{eqnarray}
denote the corresponding steady-state probability of system at state
$(n_{1},...,n_{K})$ if the network exists a stationary distribution.
In addition, let $P_{n_{iJ+k}}^{s(iJ+k)}$ denote the probability at
state $n_{iJ+k}$ for class $iJ+k$ with $i\in\{0,1,...,I-1\}$ and
$k\in\{1,2,...,J\}$.

Under the usual convention, let $x\wedge y$ denote the smaller one
of any two real numbers $x$ and $y$. Let $x\vee y$ denote the larger
one of $x$ and $y$, i.e., $x\wedge y\equiv\mbox{min}\{x,y\}$ and
$x\vee y\equiv\mbox{max}\{x,y\}$. Then, for each $k\in\{1,...,K\}$,
we have the following notation,
\begin{eqnarray}
a(x,y)\equiv\left(x\wedge(c_{s(k)}-y)\right)\vee 0.
\elabel{axy}
\end{eqnarray}
Finally, for $i>0$, define
\begin{eqnarray}
\kappa(n_{iJ+k})\equiv
   \left\{\begin{array}{ll}
            0  & \mbox{if}\;\;n_{iJ+k}=0,\\
            \sum_{\sum_{u=0}^{i-1}n_{uJ+k}<c_{s(k)}}
            a\left(n_{iJ+k},\sum_{u=0}^{i-1}n_{uJ+k}
            \right)
            \prod_{u=0}^{i-1}P_{n_{uJ+k}}^{s(uJ+k)}
            &\mbox{if}\;\;\;n_{iJ+k}>0.
          \end{array}
   \right.
\nonumber
\end{eqnarray}
Furthermore, for $i=0$, define $\kappa(n_{k})\equiv n_{k}\wedge
c_{s(k)}$.

\begin{theorem}\label{mainth}
{\bf (Steady-State Distribution)} Assume that all of the service
times and external interarrival times are exponentially distributed
with rates as before, and the traffic intensity $\rho_{j}<1$ for all
$j\in\{1,...,J\}$ . Furthermore, if
\begin{eqnarray}
\frac{\lambda_{iJ+k}} {\kappa(n_{iJ+k})\mu_{iJ+k}}<1 \elabel{kvp}
\end{eqnarray}
for $i\in\{0,...,I-1\}$, then, for each
$h\in\{1,2,...,I\}$, the
steady-state distribution is given by the
following product forms
\begin{eqnarray}
P_{n_{(1,hJ)}}=\prod_{i=0}^{h-1}
\prod_{k=1}^{J}P_{n_{iJ+k}}^{s(iJ+k)}. \elabel{product}
\end{eqnarray}
More precisely, for $n_{iJ+k}\geq 1$,
\begin{eqnarray}
P_{n_{iJ+k}}^{s(iJ+k)}=P_{0}^{s(iJ+k)}
\prod_{r_{iJ+k}=1}^{n_{iJ+k}}
\lambda_{iJ+k}/\kappa(r_{iJ+k})\mu_{iJ+k},
\elabel{pnk}
\end{eqnarray}
and for $n_{iJ+k}=0$, the initial distribution
$P^{s(iJ+k)}_{0}$ is determined by
\begin{eqnarray}
P^{s(iJ+k)}_{0}=
\frac{1}{1+\sum_{n_{iJ+k}=1}^{\infty}
\prod_{r_{iJ+k}=1}^{n_{iJ+k}}
\frac{\lambda_{iJ+k}}
{\kappa(r_{iJ+k})\mu_{iJ+k}}}.
\nonumber
\end{eqnarray}
\end{theorem}

The network stability condition \eq{kvp} is somehow complicated
although it can be checked iteratively by computer in practical
usages. The following theorems relate it to primitive interarrival
time rates and service rates.

\begin{proposition}
{\bf (Network Stability Condition)}\label{sn}
Under the exponential assumptions as stated in
Theorem~\ref{mainth}, if $\rho_{j}<1$ for each
$j\in\{1,...,J\}$, the stability condition
\eq{kvp} holds for the following networks:

{\bf Net I} Multiclass networks with single-server stations, i.e.,
the number $c_{j}$ of servers is one for all stations
$j\in\{1,...,J\}$ while the number $I$ of job types can be
arbitrary.

{\bf Net II} The number $c_{j}$ of servers can
be arbitrary for all stations $j\in\{1,...,J\}$ while
the number of job types equals two ($I=2$).
\end{proposition}

We conjecture that $\rho_{j}<1$ for $j\in\{1,...,J\}$ implies the
condition \eq{kvp} for our general network. Nevertheless, owing to
complex computation involved, the corresponding analytical
illustration is not a trivial task.

\begin{example}\label{exampleI}
Consider a network with three job types ($I=3$) and at least one
station $j=s(k)$ having two servers ($c_{s(k)}=2$ for
$k\in\{1,...,J\}$) while other stations having at most two servers
($c_{s(l)}\leq 2$ for $l\in\{1,2,...,J\}\setminus\{k\}$). For a
station $s(k)$ with three servers, the condition \eq{kvp} can be
explicitly expressed as follows,
\begin{eqnarray}
\frac{\frac{\lambda_{2J+k}}{\mu_{2J+k}}
\left(1+\frac{\lambda_{k}}{\mu_{k}}
+\frac{\lambda_{J+k}}{2\mu_{J+k}}
-\left(\frac{\lambda_{k}}{2\mu_{k}}\right)^{2} -\frac{1}{4}
\left(\frac{\lambda_{k}}{\mu_{k}}\right)^{3}
+\frac{\lambda_{k}\lambda_{J+k}}{4\mu_{k}\mu_{J+k}} \right)}
{2\left(1-\frac{\lambda_{k}}{2\mu_{k}}-
\frac{\lambda_{J+k}}{2\mu_{J+k}}\right)
\left(1+\frac{\lambda_{k}}{\mu_{k}}
+\frac{\lambda_{J+k}}{2\mu_{J+k}}
-\left(\frac{\lambda_{k}}{2\mu_{k}}\right)^{2} -\frac{1}{4}
\left(\frac{\lambda_{k}}{\mu_{k}}\right)^{3}
+2\left(1-\frac{\lambda_{k}}{2\mu_{k}}\right)
\frac{\lambda_{k}\lambda_{J+k}} {4\mu_{k}\mu_{J+k}}\right)} <1
\nonumber
\end{eqnarray}
for $n_{2J+k}\geq 2$. Under $\rho_{j}<1$, it is easy to see that the
above inequality is true if $\lambda_{k}/\mu_{k}\leq 1$. Numerical
tests~in Table~\ref{testscd} have also been conducted and show that
the inequality is true even for $1<\lambda_{k}/\mu_{k}<2$, but the
corresponding analytic demonstration could be nontrivial. The
detailed illustration of the example will be provided at the end of
this paper.
\begin{table}[tbh]
\begin{center}
\begin{tabular}{|c|p{0.5in}|p{0.5in}|p{0.5in}|p{0.5in}|p{0.5in}|p{0.5in}|}
\hline\hline \multicolumn{4}{|c|}{$1<\lambda_{k}/\mu_{k}<2$,
$\rho_{s(k)}<1$}\\ \hline\hline
\multicolumn{1}{|c|}{$\lambda_{2J+k}/\kappa(n_{2J+k})\mu_{2J+k}$~~}
& \multicolumn{1}{|c|}{$\lambda_{k}/\mu_{k}$~~~~~~~~~~~~~~} &
\multicolumn{1}{|c|}{$\lambda_{J+k}/\mu_{J+k}$~~~~~~~~~~} &
\multicolumn{1}{|c|}{$\lambda_{2J+k}/\mu_{2J+k}$~~~~~~~~} \\
\hline\hline 0.16748&$\hskip 0.3cm 1.10$&$0.30$&$0.10$  \\ \hline
0.98814&$\hskip 0.3cm 1.10$&$0.30$&$0.59$  \\ \hline 0.99111&$\hskip
0.3cm 1.50$&$0.30$&$0.19$  \\ \hline 0.85890&$\hskip 0.3cm
1.90$&$0.05$&$0.04$  \\ \hline\hline
\end{tabular}
\end{center}
\caption{Numerical tests for network stability condition \eq{kvp}}
\label{testscd}
\end{table}
\end{example}

\begin{remark}
The justifications of Theorem~\ref{mainth} and Proposition~\ref{sn}
are postponed to Section~\ref{proof}. Instead, we will first use
these results as comparison criteria to illustrate the effectiveness
of the diffusion approximation models developed in
Dai~\cite{dai:difapp} and answer the question on when these
approximation models can be employed.
\end{remark}

\begin{remark}
Under a general Whittle network framework, the exact
solutions are presented by Serfozo~\cite{ser:intsto}
for some multiclass networks, which include
those with sector-dependent and class-station-dependent
service rates such as BCMP networks introduced by Baskett,
Chandy, Muntz, and Palacios~\cite{bascha:opeclo}. Without
considering the interaction among different stations, the
distinguishing feature of these networks is that a job's
service rate at a station may consist of two intensities:
one referred as {\em station intensity} is a
function of the total queue length at the station, and the
other one referred as {\em class intensity} is a
function of the queue length in the same class as the job
being served. However, for a station in our networks, the
station intensity is not only a function of the total queue
length but also a function of combinations of queue lengths
of various classes, and the class intensity depends not
only the queue length of itself and/or the total queue
length but also the numbers of jobs in other classes at the
station.
\end{remark}

\section{Numerical comparisons}\label{comparison}

First of all, we note that partial results presented in this section
was briefly reported in the short conference version
(Dai~\cite{dai:profor}) of this paper. More precisely, we consider a
network with single-server station and under preemptive priority
service discipline. By employing the exact solutions developed in
previous sections, we conduct performance comparisons between these
product-form solutions and the approximating ones of Brownian
network models.

Brownian network models, known as semimartingale reflecting Brownian
motions (SRBM), have been widely employed as approximating models
for multiclass queueing networks with general interarrival and
service time distributions when the traffic intensity defined in
\eq{rhoj} is suitably large or close to one (see, e.g.,
Dai~\cite{dai:difapp}).

For the network with exponential interarrival and service time
distributions, by using Theorem~\ref{mainth}, we get the
steady-state mean queue length for each class $k\in\{1,...,K\}$ as
\begin{eqnarray}
EQ_{k}(\infty) =\sum_{n_{k}=0}^{\infty}n_{k}P_{n_{k}}^{s(k)}
=\frac{\lambda_{k}m_{k}}{1-\sum_{i=1}^{k} \lambda_{i}m_{i}}
=\frac{\alpha_{k}m_{k}/(1-p_{kk})}
{1-\sum_{i=1}^{k}\alpha_{i}m_{i}/(1-p_{ii})} \elabel{exactlength}
\end{eqnarray}
and the expected total time (sojourn time) a job had to spend in the
system as
\begin{eqnarray}
ET_{k}(\infty) =\frac{m_{k}}{1-\sum_{i=1}^{k} \lambda_{i}m_{i}}
=\frac{m_{k}} {1-\sum_{i=1}^{k}\alpha_{i}m_{i}/(1-p_{ii})}.
\elabel{exactdelay}
\end{eqnarray}

For the network with general interarrival and service time
distributions, owing to the nature of our network routing structure,
the higher priority classes are independent of lower ones. Then, it
follows from the studies in Dai~\cite{dai:difapp},
Harrison~\cite{har:bromot}, and Chen and Yao~\cite{cheyao:funque}
that the steady-state mean sojourn time and mean queue length for
each class $k$ can iteratively be calculated with respect to
priority order as follows,
\begin{eqnarray}
ET_{k}(\infty)=\frac{EW_{k}(\infty)+m_{k}}
{1-\sum_{i=1}^{k-1}\lambda_{i}m_{i}}\;\;\;\;\mbox{and}\;\;\;
EQ_{k}(\infty)=\lambda_{k}ET_{k}(\infty), \elabel{mequbn}
\end{eqnarray}
where $EW_{k}(\infty)$ is the steady-state mean total workload for
all classes $i\in\{1,...,k\}$ and $k=1,...,K$. More precisely, it is
given by
\begin{eqnarray}
EW_{k}(\infty)=\sigma_{k}^{2}/2|\theta_{k}|, \elabel{meanr}
\end{eqnarray}
where
\begin{eqnarray}
&&\theta_{k}=\left(1-p_{kk}\right)
\left(\sum_{i=1}^{k}\lambda_{i}m_{i}-1\right),
\elabel{thetar}\\
&&\sigma_{k}^{2}=(1-p_{kk})^{2}\sum_{i=1}^{k}
\left(\lambda_{i}b_{i}+\frac{m_{i}^{2}
\left(\alpha_{i}^{3}a_{i}+\lambda_{i}p_{ii}(1-p_{ii})\right)}
{(1-p_{ii})^{2}}\right), \elabel{sigmar}
\end{eqnarray}
and $a_{i}$, $b_{i}$ are the variances of interarrival and service
time sequences for each class $i$.

In our numerical comparisons, we consider an exponential network
with $K=3$. For this case, our corresponding data are listed in
Table~\ref{single2}. In the table, ERROR=SRBM-EXACT and
RATIO=($|$ERROR$|$/EXACT)*100$\%$,
$E\hat{Q}(\infty)=\sum_{i=1}^{3}EQ_{i}(\infty)$ and
$E\hat{T}_{k}(\infty)=\sum_{i=1}^{k}T_{i}(\infty)$ for $k=1,2,3$.
From the table, we can see that SRBMs are more reasonable
approximations of their physical queueing counterparts when traffic
intensities for the lowest priority jobs are relatively large.

\begin{table}
\begin{center}
\begin{tabular}{|c|p{0.5in}|p{0.5in}|p{0.5in}|p{0.5in}|p{0.5in}|p{0.5in}|p{0.5in}|p{0.5in}|}
\hline\hline \multicolumn{9}{|c|}{$K=3$, $m_{1}=1$, $m_{2}=2$,
$m_{3}=3$,
                     $p_{11}=0.1$, $p_{22}=0.2$, $p_{33}=0.3$}\\
\hline\hline\hline \multicolumn{1}{|c|}{~~~~~~}        &
\multicolumn{1}{|c|}{$\lambda_{1}$} &
\multicolumn{1}{|c|}{$\lambda_{2}$} &
\multicolumn{1}{|c|}{$\lambda_{3}$} & \multicolumn{1}{|c|}{$\rho$}
& \multicolumn{1}{|c|}{$E\hat{Q}(\infty)$}     &
\multicolumn{1}{|c|}{$E\hat{T}_{1}(\infty)$} &
\multicolumn{1}{|c|}{$E\hat{T}_{2}(\infty)$} &
\multicolumn{1}{|c|}{$E\hat{T}_{3}(\infty)$} \\ \hline\hline\hline
EXACT&$0.20$&$0.20$&$0.06$&$0.78$&$ 2.0682$&$1.2500$&$ 6.2500$&$
14.8864$\\ \hline SRBM &$0.20$&$0.20$&$0.06$&$0.78$&$
2.2410$&$1.2500$&$ 6.4722$&$ 17.0253$\\ \hline ERROR&$    $&$    $&$
$&$    $&$-0.1728$&$0.0000$&$ 0.2222$&$ 2.1389 $\\ \hline RATIO&$
$&$    $&$    $&$    $&$ 8.35\%$&$0.00\%$&$ 3.56\%$&$ 14.37\%$\\
\hline\hline\hline
EXACT&$0.20$&$0.20$&$0.085$&$0.855$&$ 3.0086$&$1.2500$&$ 6.2500$&$
21.9397$\\ \hline SRBM &$0.20$&$0.20$&$0.085$&$0.855$&$
3.2349$&$1.2500$&$ 6.4722$&$ 24.0785$\\ \hline ERROR&$    $&$    $&$
$&$     $&$-0.2263$&$0.0000$&$ 0.2222$&$ 2.1389 $\\ \hline RATIO&$
$&$    $&$     $&$     $&$ 7.52\%$&$0.00\%$&$ 3.56\%$&$ 9.75\% $\\
\hline\hline\hline
EXACT&$0.20$&$0.20$&$0.115$&$0.945$&$ 7.5227$&$1.2500$&$ 6.2500$&$
55.7955$\\ \hline SRBM &$0.20$&$0.20$&$0.115$&$0.945$&$
7.8131$&$1.2500$&$ 6.4722$&$ 57.9344$\\ \hline ERROR&$    $&$    $&$
$&$     $&$-0.2904$&$0.0000$&$ 0.2222$&$ 2.1389 $\\ \hline RATIO&$
$&$    $&$     $&$     $&$ 3.86\%$&$0.00\%$&$ 3.56\%$&$ 3.83\% $\\
\hline\hline\hline
EXACT&$0.10$&$0.15$&$0.115$&$0.745$&$ 1.9641$&$1.1111$&$ 4.4444$&$
12.8758$\\ \hline SRBM &$0.10$&$0.15$&$0.115$&$0.745$&$
2.0944$&$1.1111$&$ 4.5432$&$ 13.8804$\\ \hline ERROR&$    $&$    $&$
$&$     $&$-0.1303$&$0.0000$&$ 0.0988$&$ 1.0046 $\\ \hline RATIO&$
$&$    $&$     $&$     $&$ 6.64\%$&$0.00\%$&$ 2.22\%$&$ 7.80\% $\\
\hline\hline\hline
EXACT&$0.10$&$0.15$&$0.15$&$0.85$&$ 3.6111$&$1.2500$&$ 4.4444$&$
21.1111$\\ \hline SRBM &$0.10$&$0.15$&$0.15$&$0.85$&$
3.7766$&$1.2500$&$ 4.5432$&$ 22.1157$\\ \hline ERROR&$    $&$    $&$
$&$    $&$-0.1655$&$0.0000$&$ 0.0988$&$ 1.0046 $\\ \hline RATIO&$
$&$    $&$    $&$    $&$ 4.58\%$&$0.00\%$&$ 2.22\%$&$ 4.76\% $\\
\hline\hline\hline
EXACT&$0.10$&$0.15$&$0.185$&$0.955$&$ 12.9445$&$1.2500$&$ 4.4444$&$
67.7778$\\ \hline SRBM &$0.10$&$0.15$&$0.185$&$0.955$&$
13.1451$&$1.2500$&$ 4.5432$&$ 68.7824$\\ \hline ERROR&$    $&$
$&$     $&$     $&$-0.2007 $&$0.0000$&$ 0.0988$&$ 1.0046 $\\ \hline
RATIO&$    $&$    $&$     $&$     $&$ 1.55\% $&$0.00\%$&$ 2.22\%$&$
1.48\% $\\ \hline\hline\hline
EXACT&$0.45$&$0.22$&$0.02$&$0.95$&$ 6.0182$&$1.8182$&$ 20.0000$&$
61.8182$\\ \hline SRBM &$0.45$&$0.22$&$0.02$&$0.95$&$
6.3818$&$1.8182$&$ 20.7273$&$ 72.0000$\\ \hline ERROR&$    $&$
$&$    $&$    $&$-0.3636$&$0.0000$&$ 0.7273 $&$ 10.1818$\\ \hline
RATIO&$    $&$    $&$    $&$    $&$ 6.04\%$&$0.00\%$&$ 3.64\% $&$
16.47\% $\\ \hline\hline\hline
EXACT&$0.45$&$0.22$&$0.03$&$0.98$&$ 9.3182$&$1.8182$&$ 20.0000$&$
151.8181$\\ \hline SRBM &$0.45$&$0.22$&$0.03$&$0.98$&$
9.7836$&$1.8182$&$ 20.7273$&$ 161.9999$\\ \hline ERROR&$    $&$
$&$    $&$    $&$-0.4655$&$0.0000$&$ 0.7273 $&$ 10.1818 $\\ \hline
RATIO&$    $&$    $&$    $&$    $&$ 5.00\%$&$0.00\%$&$ 3.64\% $&$
6.71\%  $\\ \hline\hline
\end{tabular}
\end{center}
\caption{Performance comparisons for a priority multiclass network
with three job types} \label{single2}
\end{table}

\section{Proofs of Theorem \ref{mainth}
and Proposition~\ref{sn}}
\label{proof}

\subsection{Proof of Theorem \ref{mainth}}

For convenience, we introduce some additional notations. Let $\Delta
E_{k}(t)$, $\Delta S_{k}(t)$ and $\Delta A_{k}(t)$ be defined by
\begin{eqnarray}
&&\Delta E_{k}(t)\equiv E_{k}(t+\Delta t)-E_{k}(t),
\elabel{deltaE}\\
&&\Delta S_{k}(t)\equiv S_{k}(t+\Delta t)-S_{k}(t),
\elabel{deltaS}\\
&&\Delta A_{k}(t)\equiv A_{k}(t+\Delta t)-A_{k}(t),
\elabel{deltaA}
\end{eqnarray}
which denote the cumulative external arrivals to class $k$,
the cumulative number of jobs finished services at class $k$,
and the total arrivals to class $k$ in $[t,t+\Delta t]$.
Then, we can justify Theorem \ref{mainth} by induction
as in the following several steps.

{\bf Step One.} We consider the steady-state distribution for the
highest rank classes with each index $k\in\{1,2,...,J\}$. In this
case, the type index $i=0$ as used in Theorem~\ref{mainth}. Owing to
the preemptive-resume service discipline and the class routing
structure, we know that these $J$ classes form a Jackson network.
Then, by the theorem in Jackson \cite{jac:netwai}, we have the
following product form
\begin{eqnarray}
P_{n_{(1,J)}}=P^{s(1)}_{n_{1}}P^{s(2)}_{n_{2}} \cdot\cdot\cdot
P^{s(J)}_{n_{J}}, \elabel{highestD}
\end{eqnarray}
where $P_{n_{k}}^{s(k)}$ denotes the steady-state probability at
state $n_{k}$ for class $k\;(\in\{1,2,...,J\})$ at station $s(k)$.
More precisely, it is given by
\begin{eqnarray}
P_{n_{k}}^{s(k)}=\left\{\begin{array}{ll}
                 P_{0}^{s(k)}(\lambda_{k}/\mu_{k})^{n_{k}}/n_{k}!\;,
                   & \;\;\mbox{if}\;\;0\leq n_{k}<c_{k},\\
                 P_{0}^{s(k)}(\lambda_{k}/\mu_{k})^{n_{k}}
                 /n_{k}!(n_{k})^{n_{k}-c_{k}},
                   & \;\;\mbox{if}\;\;n_{k}\geq c_{k}
              \end{array}
         \right.
\elabel{hI}
\end{eqnarray}
with $P_{0}^{s(k)}$ being determined by equation
$\sum_{n_{k}=0}^{\infty}P_{n_{k}}^{s(k)}=1$, i.e.,
\begin{eqnarray}
P_{0}^{s(k)}=\frac{1}{\sum_{n_{k}=0}^{c_{s(k)}-1}
\frac{1}{n_{k}!}
\left(\frac{\lambda_{k}}{\mu_{k}}\right)^{n_{k}}
+\frac{\left(\lambda_{k}/\mu_{k}\right)^{c_{s(k)}}}
{c_{s(k)}!\left(1-\lambda_{k}/c_{s(k)}\mu_{k}\right)}}.
\elabel{hII}
\end{eqnarray}

{\bf Step Two.} We derive the steady-state distribution for the
highest rank and the second highest rank classes with
$k\in\{1,2,...,2J\}$. Owing to the preemptive-resume service
discipline and the job routing structure, we have
\begin{eqnarray}
&&P_{n_{(1,2J)}}(t)
\elabel{P(t)2J}\\
&&=P_{n_{(1,J)}/n_{(J+1,2J)}}(t)P_{(n_{(J+1,2J)}}(t)
\nonumber\\
&&=P_{n_{(1,J)}}(t)P_{n_{(J+1,2J)}}(t),
\nonumber
\end{eqnarray}
where $P_{n_{(1,J)}/n_{(J+1,2J)}}(t)$ is the conditional probability
at time $t$ for classes with $k\in\{1,...,J\}$ at state $n_{(1,J)}$
in terms of classes with $k\in\{J+1,...,2J\}$ at state
$n_{(J+1,2J)}$, i.e.,
\begin{eqnarray}
P_{n_{(1,J)}/n_{(J+1,2J)}}(t)
=P\left\{Q_{(1,J)}(t)=
n_{(1,J)}/(Q_{(J+1,2J)}(t)=n_{(J+1,2J)}\right\}.
\nonumber
\end{eqnarray}
In order to get the steady state distribution for
$P_{n_{(J+1,2J))}}(t)$, we consider each state
$n_{(J+1,2J)}$ at time $t$ for the second highest rank
class jobs. There are several ways in which the system
can reach it. They can be summarized in the following
formula,
\begin{eqnarray}
&&P_{n_{(J+1,2J)}}(t+\Delta t)
\elabel{pdeltat}\\
&&=\left(1-\sum_{k=1}^{J}\alpha_{J+k}
\Delta t-\sum_{k=1}^{J}
\sum_{n_{1},...,n_{J}=0}^{\infty}
a(n_{J+k},n_{k})\mu_{J+k}(1-p_{J+k,J+k})
P_{n_{(1,J)}}(t)\Delta t\right)
\nonumber\\
&&\;\;\;\;\;\;\;\;P_{n_{(J+1,2J)}}(t)
\nonumber\\
&&\;\;\;\;+\sum_{k=1}^{J}
\left(\sum_{n_{1},...,n_{J}=0}^{\infty}
a(n_{J+k}+1,n_{k})\mu_{J+k}
\left(1-\sum_{l=1}^{J}p_{J+k,J+l}\right)
P_{n_{(1,J)}}(t)\Delta t\right)
\nonumber\\
&&\;\;\;\;\;\;\;\;\;\;\;\;\;\;\;
P_{n_{J+1},...,n_{J+k}+1,...,n_{2J}}(t)
\nonumber\\
&&\;\;\;\;+\sum_{k=1}^{J}
(n_{J+k}\wedge 1)\alpha_{J+k}
\Delta t P_{n_{J+1},...,n_{J+k}-1,...,n_{2J}}(t)
\nonumber\\
&&\;\;\;\;+\sum_{r,s=1,s\neq r}^{J}
\left(\sum_{n_{1},...,n_{J}=0}^{\infty}
a(n_{J+r}+1,n_{r})
\mu_{J+r}p_{J+r,J+s}P_{n_{(1,J)}}(t)\right)
\nonumber\\
&&\;\;\;\;\;\;\;\;\;\;\;\;\;\;\;\;\;\;\;\;\;
(n_{J+s}\wedge 1)\Delta t
P_{n_{J+1},...,n_{J+r}+1,...,n_{J+s}-1,...,n_{2J}}(t)
\nonumber\\
&&\;\;\;\;+o(\Delta t),
\nonumber
\end{eqnarray}
where $a(n_{J+k},n_{k})=(n_{J+k} \wedge(c_{s(k)}-n_{k}))\vee 0$ as
defined in \eq{axy}. $o(\Delta t)$ is infinitesimal in terms of
$\Delta t$, i.e., $o(\Delta t)/\Delta t\rightarrow 0$ as $\Delta
t\rightarrow\infty$. The equation \eq{pdeltat} can be illustrated in
the following disjoint events.

{\bf Event ${\cal A}$:} The $J$-dimensional queue length process
$Q_{(J+1,2J)}(t)$ keeps at state $n_{(J+1,2J)}$ unchanging at times
$t$ and $t+\Delta t$, i.e.,
$$Q_{(J+1,2J)}(t)=n_{(J+1,2J)}
,\;\;\;Q_{(J+1,2J)}(t+\Delta t)=n_{(J+1,2J)}.$$
This event involves the following two parts:

{\em Part One.} $\Delta S_{J+k}(t)=\Delta A_{J+k}(t)=0$ for all
$k\in\{1,2,...,J\}$, i.e., no external jobs arrive to classes with
indices belonging to $\{J+1,...,2J\}$ during $[t,t+\Delta t]$ while
no jobs finish their services, either because the jobs being served
at time $t$ require longer service times than $\Delta t$ or because
the services are blocked or interrupted by higher rank class jobs.
Then, the probability for Part One can be expressed in the following
product form,
\begin{eqnarray}
&&P\left\{\sum_{k=1}^{J}\Delta S_{J+k}(t)
=\sum_{k=1}^{J}\Delta A_{J+k}(t)=0\right\}
\elabel{partone}\\
&&=P\left\{\sum_{k=1}^{J}\Delta S_{J+k}(t)
=\sum_{k=1}^{J}\Delta E_{J+k}(t)=0\right\}
\nonumber\\
&&=\prod_{k=1}^{J}P\left\{\Delta S_{J+k}(t)=
\Delta E_{J+k}(t)=0/\sum_{s=k}^{J}
\Delta S_{J+s+1}(t)
=\sum_{s=k}^{J}\Delta E_{J+s+1}(t)=0\right\}
\nonumber\\
&&=\prod_{k=1}^{J}P\left\{{\cal E}_{J+k}(t)\right\}
\nonumber\\
&&=\prod_{k=1}^{J}\left(1-\alpha_{J+k}
\Delta t-a_{J+k}(t)\mu_{J+k}\Delta t
+ o(\Delta t)\right)
\nonumber\\
&&=1-\sum_{k=1}^{J}\alpha_{J+k}\Delta t-
\sum_{k=1}^{J}a_{J+k}(t)\mu_{J+k}\Delta t
+o(\Delta t),
\nonumber
\end{eqnarray}
where $a_{J+k}(t)$ is defined to be
\begin{eqnarray}
a_{J+k}(t)\equiv\left(Q_{J+k}(t)
\wedge(c_{s(k)}-Q_{k}(t))\right)\vee 0.
\elabel{ajkt}
\end{eqnarray}
The event ${\cal E}_{J+k}(t)$ in the third equality
of \eq{partone} is defined as follows,
\begin{eqnarray}
{\cal E}_{J+k}(t)=\left\{\Delta S_{J+k}(t)=
\Delta E_{J+k}(t)=0/\sum_{s=k}^{J}
\Delta S_{J+s+1}(t)
=\sum_{s=k}^{J}\Delta E_{J+s+1}(t)=0\right\}.
\nonumber
\end{eqnarray}
To explain the fourth equality of \eq{partone}, we introduce more
notations. Let $b_{k}(t)\Delta t$ denote the probability that
$\Delta A_{k}(t)=1$ for $k\in\{1,...,J\}$ , and $c_{J+k}(t)\Delta t$
be the probability that $\sum_{s=1}^{k-1}\Delta S_{J+s}(t)=1$ for
$k\in\{2,3,...,J\}$, i.e.,
\begin{eqnarray}
&&b_{k}(t)\Delta t\equiv P\left\{\Delta A_{k}(t)=1\right\},
\nonumber\\
&&c_{J+k}(t)\Delta t\equiv P\left\{\sum_{s=1}^{k-1}
\Delta S_{J+s}(t)=1\right\}.
\nonumber
\end{eqnarray}
These probabilities can be explicitly expressed in terms of
the external arrival rates, service rates and network states,
for example,
\begin{eqnarray}
b_{k}(t)\Delta t=\left(\alpha_{k}+
\sum_{l\in\{1,...,J\}}p_{lk}\mu_{l}
\left(c_{s(l)}\wedge Q_{l}(t)\right)\right)\Delta t.
\nonumber
\end{eqnarray}
Thus, by the independent assumptions on external
arrival and service processes among different stations and
classes, and for each $k\in\{2,3,...,J\}$, we have
\begin{eqnarray}
&&P\left\{{\cal E}_{J+k}(t)\right\}
\elabel{ejk}\\
&&=\sum_{n=0}^{\infty}\sum_{m=0}^{\infty}
P\left\{{\cal E}_{J+k}(t)/
\Delta A_{k}(t)=n,\sum_{s=1}^{k-1}
\Delta S_{J+s}(t)=m\right\}
\nonumber\\
&&\;\;\;\;\;\;\;\;\;\;\;\;\;\;\;\;\;
P\left\{\Delta A_{k}(t)=n\right\}
P\left\{\sum_{s=1}^{k-1}\Delta S_{J+s}(t)=n\right\}
\nonumber\\
&&=P\left\{{\cal E}_{J+k}(t)/
\Delta A_{k}(t)=0,\sum_{s=1}^{k-1}
\Delta S_{J+s}(t)=0\right\}
P\left\{\Delta A_{k}(t)=0\right\}
P\left\{\sum_{s=1}^{k-1}\Delta S_{J+s}(t)=0\right\}
\nonumber\\
&&\;+P\left\{{\cal E}_{J+k}(t)/
\Delta A_{k}(t)=1,\sum_{s=1}^{k-1}
\Delta S_{J+s}(t)=0\right\}
P\left\{\Delta A_{k}(t)=1\right\}
P\left\{\sum_{s=1}^{k-1}\Delta S_{J+s}(t)=0\right\}
\nonumber\\
&&\;+P\left\{{\cal E}_{J+k}(t)/
\Delta A_{k}(t)=0,\sum_{s=1}^{k-1}
\Delta S_{J+s}(t)=1\right\}
P\left\{\Delta A_{k}(t)=0\right\}
P\left\{\sum_{s=0}^{k-1}\Delta S_{J+s}(t)=1\right\}
\nonumber\\
&&\;+o(\Delta t)
\nonumber\\
&&=\left(1-a_{J+k}(t)\mu_{J+k}\Delta t\right)
\left(1-\alpha_{J+k}\Delta t\right)
\left(1-b_{k}(t)\Delta t\right)
\left(1-c_{J+k}(t)\Delta t\right)
\nonumber\\
&&\;\;\;+\left(1-\left(\left(Q_{J+k}(t)\wedge
\left(c_{s(k)}-(Q_{k}(t)+1)\right)\right)
\vee 0\right)\mu_{J+k}\Delta t\right)
\left(1-\alpha_{J+k}\Delta t\right)
\nonumber\\
&&\;\;\;\;\;\;\left(b_{k}(t)\Delta t-p_{kk}\mu_{k}
\left(c_{s(k)}\wedge Q_{k}(t)\right)\Delta t\right)
\left(1-c_{J+k}(t)\Delta t\right)
\nonumber\\
&&\;\;\;+\left(1-\left(\left(Q_{J+k}(t)\wedge
(c_{s(k)}-Q_{k}(t))\right)\vee 0\right)
\mu_{J+k}\Delta t\right)\left(1-\alpha_{J+k}
\Delta t\right)
\nonumber\\
&&\;\;\;\;\;\;p_{kk}\mu_{k}\left(c_{s(k)}\wedge
Q_{k}(t)\right)
\Delta t\left(1-c_{J+k}(t)\Delta t\right)
\nonumber\\
&&\;\;\;+\left(1-\left(\left(\left(Q_{J+k}(t)+1
\right)\wedge\left(c_{s(k)}-Q_{k}(t)\right)
\right)\vee 0\right)\mu_{J+k}\Delta t\right)
\left(1-\alpha_{J+k}\Delta t\right)
\nonumber\\
&&\;\;\;\;\;\;\left(1-b_{k}(t)\Delta t\right)
c_{J+k}(t)\Delta t
\nonumber\\
&&\;\;\;+o(\Delta t)
\nonumber\\
&&=1-\alpha_{J+k}\Delta t-a_{J+k}(t)
\mu_{J+k}\Delta t + o(\Delta t).
\nonumber
\end{eqnarray}
For the case $k=1$, one can use the similar way
to check that the result in the above equation
is also true.

{\em Part Two.} $\Delta S_{J+k}(t)=\Delta A_{J+k}(t)\geq 1$ for at
least one $k\in\{1,2,...J\}$, i.e., the number (at least one) of
jobs finished their services for class $J+k$ in $[t,t+\Delta t]$
equals that of jobs arrived at the class. It is easy to see that the
probability corresponding part one is given by
\begin{eqnarray}
&&P\left\{\Delta S_{J+k}(t)=\Delta A_{J+k}(t)\geq 1\;\;
\mbox{for at least one}\;\;k\in\{1,2,...J\}\right\}
\elabel{partwo}\\
&&=\sum_{k=1}^{J}a_{J+k}(t)\mu_{J+k}p_{J+k,J+k}\Delta t +o(\Delta
t), \nonumber
\end{eqnarray}
where $a_{J+k}(t)$ is defined in \eq{ajkt}.

Then, it follows from equations \eq{partone} and \eq{partwo} that
the probability for Event $A$ is given by
\begin{eqnarray}
&&P\left\{Q_{(J+1,2J)}(t)=n_{(J+1,2J)}
;\;Q_{(J+1,2J)}(t+\Delta t)=n_{(J+1,2J)}\right\}
\nonumber\\
&&=\sum_{n_{1},...,n_{J}=0}^{\infty}
P\left\{Q_{(J+1,2J)}(t+\Delta t)=n_{(J+1,2J)}
/Q_{(1,2J)}(t)=n_{(1,2J)}\right\}
\nonumber\\
&&\;\;\;\;\;\;\;\;\;\;\;\;\;\;\;\;\;\;\;
P\left\{Q_{(1,J)}(t)=n_{(1,J)}\right\}
P\left\{Q_{(J+1,2J)}(t)=n_{(J+1,2J)}\right\}
+o(\Delta t)
\nonumber\\
&&=\left(1-\sum_{k=1}^{J}\alpha_{J+k}
\Delta t-\sum_{k=1}^{J}
\sum_{n_{1},...,n_{J}=0}^{\infty}
a(n_{J+k},n_{k})\mu_{J+k}(1-p_{J+k,J+k})
P_{n_{(1,J)}}(t)\Delta t\right)
\nonumber\\
&&\;\;\;\;\;P_{n_{(J+1,2J)}}(t)+o(\Delta t).
\nonumber
\end{eqnarray}

{\bf Event ${\cal B}$:} There is a
$k\in\{1,...,J\}$ such that $Q_{J+k}(t)=n_{J+1}+1$,
$Q_{J+k}(t+\Delta t)=n_{J+1}$
and for all $l\in\{1,...,J\}\setminus\{k\}$,
$Q_{J+l}(t)=Q_{J+l}(t+\Delta t)=n_{J+l}$. Similar to
the discussion in Event ${\cal A}$, we can obtain
the probability for Event ${\cal B}$ as follows,
\begin{eqnarray}
&&\sum_{k=1}^{J}P\{Q_{J+k}(t)=n_{J+k}+1,
Q_{J+k}(t+\Delta t)=n_{J+k},
\nonumber\\
&&\;\;\;\;\;\;\;\;\;\;\;
Q_{J+l}(t)=Q_{J+l}(t+\Delta t)=n_{J+l},\;
l\in\{1,...,J\}\setminus\{k\}\}
\nonumber\\
&&=\sum_{k=1}^{J}P\{\Delta S_{J+k}(t)=1,
\Delta A_{J+k}(t)=0,\sum_{l=1,l\neq k}^{J}
\Delta S_{J+l}(t)=\sum_{l=1,l\neq k}^{J}
\Delta A_{J+l}(t)=0,
\nonumber\\
&&\;\;\;\;\;\;\;\;\;\;\;\;\;\;
Q_{J+k}(t)=n_{J+k}+1,
Q_{J+l}(t)=n_{J+l}\;\;\mbox{for all}\;\;
l\in\{1,...,J\}\setminus\{k\}\}
+o(\Delta t)
\nonumber\\
&&=\sum_{k=1}^{J}P\{\sum_{l=1,l\neq k}^{J}
\Delta S_{J+l}(t)=\Delta A_{J+l}(t)=0,
\Delta A_{J+k}(t)=0,Q_{J+k}(t)=n_{J+k}+1,
\nonumber\\
&&\;\;\;\;\;\;\;\;\;\;\;\;\;
Q_{J+l}(t)=n_{J+l},\;l\in\{1,...,J\}
\setminus\{k\}/\Delta S_{J+k}(t)=1\}
P\left\{\Delta S_{J+k}(t)=1\right\}+o(\Delta t)
\nonumber\\
&&=\sum_{k=1}^{J}
\left(\sum_{n_{1},...,n_{J}=0}^{\infty}
a(n_{J+k}+1,n_{k})\mu_{J+k}
\left(1-\sum_{l=1}^{J}p_{J+k,J+l}\right)
P_{n_{(1,J)}}(t)\Delta t\right)
\nonumber\\
&&\;\;\;\;\;\;\;\;\;\;
P_{n_{J+1},...,n_{J+k}+1,...,n_{2J}}(t)
+o(\Delta t).
\nonumber
\end{eqnarray}

{\bf Event ${\cal C}$:} There is a $k\in\{1,...,J\}$
such that $Q_{J+k}(t)=n_{J+k}-1$,
$Q_{J+k}(t+\Delta t)=n_{J+k}$ and for all
$l\in\{1,...,J\}\setminus\{k\}$,
$Q_{J+l}(t)=Q_{J+l}(t+\Delta t)=n_{J+l}$. Then, we
can get the probability for Event ${\cal C}$ as
follows,
\begin{eqnarray}
&&\sum_{k=1}^{J}P\{Q_{J+k}(t)=n_{J+k}-1,
Q_{J+k}(t+\Delta t)=n_{J+k},
\nonumber\\
&&\;\;\;\;\;\;\;\;\;\;\;
Q_{J+l}(t)=Q_{J+l}(t+\Delta t)=n_{J+l},\;
l\in\{1,...,J\}\setminus\{k\}\}
\nonumber\\
&&=\sum_{k=1}^{J}P\{\Delta S_{J+k}(t)=0,
\Delta E_{J+k}(t)=1,\sum_{l=1,l\neq k}^{J}
\Delta S_{J+l}(t)=\sum_{l=1,l\neq k}^{J}
\Delta A_{J+l}(t)=0,
\nonumber\\
&&\;\;\;\;\;\;\;\;\;\;\;\;\;\;
Q_{J+k}(t)=n_{J+k}-1,
Q_{J+l}(t)=n_{J+l}\;\;\mbox{for all}\;\;
l\in\{1,...,J\}\setminus\{k\}\}
+o(\Delta t)
\nonumber\\
&&=\sum_{k=1}^{J}(n_{J+k}\wedge 1)\alpha_{J+k}
\Delta t P_{n_{J+1},...,n_{J+k}-1,...,n_{2J}}(t)
+o(\Delta t).
\nonumber
\end{eqnarray}

{\bf Event ${\cal D}$:} There exist $r,s\in\{1,...,J\}$ such that
$Q_{J+r}(t)=n_{J+r}+1$, $Q_{J+r}(t+\Delta t)=n_{J+s}$,
$Q_{J+s}(t)=n_{J+s}-1$, $Q_{J+s}(t+\Delta t)=n_{J+s}$, and for all
$l\in\{1,...,J\}\setminus\{r,s\}$, $Q_{J+l}(t)=Q_{J+l}(t+\Delta
t)=n_{J+l}$. Then, we can get the probability for Event ${\cal D}$
as follows,
\begin{eqnarray}
&&\sum_{r,s=1,s\neq r}^{J}P\{Q_{J+r}(t)=n_{J+r}+1,
Q_{J+r}(t+\Delta t)=n_{J+r}; Q_{J+s}(t)=n_{J+s}-1,
\nonumber\\
&&\;\;\;\;\;\;\;\;\;\;\;\;\;\;\;\;\;\;
Q_{J+s}(t+\Delta t)=n_{J+s};
Q_{J+l}(t)=Q_{J+l}(t+\Delta t)=n_{J+l}\;
\nonumber\\
&&\;\;\;\;\;\;\;\;\;\;\;\;\;\;\;\;\;\;
\mbox{for all}\;l\in\{1,...,J\}\setminus\{r,s\}\}
\nonumber\\
&&=\sum_{r,s=1,s\neq r}^{J}P\{\Delta S_{J+r}(t)=1,
\Delta A_{J+r}(t)=0;
\Delta S_{J+s}(t)=0, \Delta A_{J+s}(t)=1;
\nonumber\\
&&\;\;\;\;\;\;\;\;\;\;\;\;\;\;\;\;\;\;\;\;\;\;
\sum_{l=1,l\neq r,s}^{J}\Delta S_{J+l}(t)=
\sum_{l=1,l\neq r,s}^{J}\Delta A_{J+l}(t)=0;
\nonumber\\
&&\;\;\;\;\;\;\;\;\;\;\;\;\;\;\;\;\;\;\;\;\;\;
Q_{J+r}(t)=n_{J+r}+1,
Q_{J+s}(t)=n_{J+s}-1,
Q_{J+l}(t)=n_{J+l}
\nonumber\\
&&\;\;\;\;\;\;\;\;\;\;\;\;\;\;\;\;\;\;\;\;\;\;
\mbox{for all}\;l\in\{1,...,J\}\setminus\{r,s\}\}
+o(\Delta t)
\nonumber\\
&&=\sum_{r,s=1,s\neq r}^{J}
\left(\sum_{n_{1},...,n_{J}=0}^{\infty}
a(n_{J+r}+1,n_{r})
\mu_{J+r}p_{J+r,J+s}P_{n_{(1,J)}}(t)\right)
\nonumber\\
&&\;\;\;\;\;\;\;\;\;\;\;\;\;\;\;\;\;\;
\left(n_{J+s}\wedge 1\right)
P_{n_{J+1},...,n_{J+r}+1,...,n_{J+s}-1,...,n_{2J}}(t)
+o(\Delta t).
\nonumber
\end{eqnarray}

\indent Now we go back to discuss equation \eq{pdeltat}. After
transferring the $P_{n_{(J+1,2J)}}(t)$ from right to left, dividing
by $\Delta t$, and taking the limit as $\Delta t$ tends to zero, we
get the following differential equations.
\begin{eqnarray}
&&\dot{P}_{n_{(J+1,2J)}}(t)
\elabel{i2de}\\
&&=\left(-\sum_{k=1}^{J}\alpha_{J+k}
-\sum_{k=1}^{J}
\sum_{n_{1},...,n_{J}=0}^{\infty}
a(n_{J+k},n_{k})\mu_{J+k}\left(1-p_{J+k,J+k}\right)
P_{n_{(1,J)}}(t)\right)
\nonumber\\
&&\;\;\;\;\;\;\;\;P_{n_{(J+1,2J)}}(t)
\nonumber\\
&&\;\;\;\;+\sum_{k=1}^{J}
\left(\sum_{n_{1},...,n_{J}=0}^{\infty}
a(n_{J+k}+1,n_{k})\mu_{J+k}
\left(1-\sum_{l=1}^{J}p_{J+k,J+l}\right)
P_{n_{(1,J)}}(t)\right)
\nonumber\\
&&\;\;\;\;\;\;\;\;\;\;\;\;\;\;\;
P_{n_{J+1},...,n_{J+k}+1,...,n_{2J}}(t)
\nonumber\\
&&\;\;\;\;+\sum_{k=1}^{J}
(n_{J+k}\wedge 1)\alpha_{J+k}
P_{n_{J+1},...,n_{J+k}-1,...,n_{2J}}(t)
\nonumber\\
&&\;\;\;\;+\sum_{r,s=1,s\neq r}^{J}
\left(\sum_{n_{1},...,n_{J}=0}^{\infty}
a(n_{J+r}+1,n_{r})
\mu_{J+r}p_{J+r,J+s}P_{n_{(1,J)}}(t)\right)
\nonumber\\
&&\;\;\;\;\;\;\;\;\;\;\;\;\;\;\;\;\;\;\;\;\;
(n_{J+s}\wedge 1)
P_{n_{J+1},...,n_{J+r}+1,...,n_{J+s}-1,...,n_{2J}}(t)
\nonumber
\end{eqnarray}
Next we show that the given distribution corresponding $h=2$ in
Theorem~\ref{mainth} is a steady-state solution of those equations
described by \eq{i2de}. It is enough to demonstrate that the
derivatives in the above equations are all made zero by setting
$P_{n_{(1,2J)}}(t)=P_{n_{(1,2J)}}$, i.e., to prove that
\begin{eqnarray}
&&\left(\sum_{k=1}^{J}\alpha_{J+k}
+\sum_{k=1}^{J}
\sum_{n_{1},...,n_{J}=0}^{\infty}
a(n_{J+k},n_{k})P_{n_{(1,J)}}\mu_{J+k}\right)
P_{n_{(J+1,2J)}}
\elabel{balanceq}\\
&&=\sum_{k=1}^{J}
\left(\sum_{n_{1},...,n_{J}=0}^{\infty}
a(n_{J+k}+1,n_{k})P_{n_{(1,J)}}\mu_{J+k}
\left(1-\sum_{l=1}^{J}p_{J+k,J+l}\right)\right)
\nonumber\\
&&\;\;\;\;\;\;\;\;\;\;\;\;\;\;\;
P_{n_{J+1},...,n_{J+k}+1,...,n_{2J}}
\nonumber\\
&&\;\;\;\;+\sum_{k=1}^{J}
(n_{J+k}\wedge 1)\alpha_{J+k}
P_{n_{J+1},...,n_{J+k}-1,...,n_{2J}}
\nonumber\\
&&\;\;\;\;+\sum_{r,s=1,s\neq r}^{J}
\left(\sum_{n_{1},...,n_{J}=0}^{\infty}
a(n_{J+r}+1,n_{r})P_{n_{(1,J)}}
\mu_{J+r}p_{J+r,J+s}\right)
\nonumber\\
&&\;\;\;\;\;\;\;\;\;\;\;\;\;\;\;\;\;\;\;\;\;
(n_{J+s}\wedge 1)
P_{n_{J+1},...,n_{J+r}+1,...,n_{J+s}-1,...,n_{2J}}
\nonumber\\
&&\;\;\;\;+\sum_{k=1}^{J}
\sum_{n_{1},...,n_{J}=0}^{\infty}
a(n_{J+k},n_{k})P_{n_{(1,J)}}\mu_{J+k}p_{J+k,J+k}
P_{n_{(J+1,2J)}}.
\nonumber
\end{eqnarray}
By the definition of $\kappa(\cdot)$ in Theorem~\ref{mainth}, i.e.,
for each $k\in\{1,...,J\}$,
\begin{eqnarray}
\kappa(n_{J+k})=\sum_{n_{1},...,n_{J}=0}^{\infty}
a(n_{J+k},n_{k})P_{n_{(1,J)}},
\nonumber
\end{eqnarray}
then, we can rewrite \eq{balanceq} as follows,
\begin{eqnarray}
&&\left(\sum_{k=1}^{J}\alpha_{J+k}
+\sum_{k=1}^{J}\kappa(n_{J+k})\mu_{J+k}\right)
P_{n_{(J+1,2J)}}
\elabel{balanck}\\
&&=\sum_{k=1}^{J}\kappa(n_{J+k}+1)\mu_{J+k}
\left(1-\sum_{l=1}^{J}p_{J+k,J+l}\right)
P_{n_{J+1},...,n_{J+k}+1,...,n_{2J}}
\nonumber\\
&&\;\;\;\;+\sum_{k=1}^{J}
(n_{J+k}\wedge 1)\alpha_{J+k}
P_{n_{J+1},...,n_{J+k}-1,...,n_{2J}}
\nonumber\\
&&\;\;\;\;+\sum_{r,s=1,s\neq r}^{J}
\kappa(n_{J+r}+1)\mu_{J+r}p_{J+r,J+s}(n_{J+s}\wedge 1)
P_{n_{J+1},...,n_{J+r}+1,...,n_{J+s}-1,...,n_{2J}}
\nonumber\\
&&\;\;\;\;+\sum_{k=1}^{J}\kappa(n_{J+k})
\mu_{J+k}p_{J+k,J+k}P_{n_{(J+1,2J)}}.
\nonumber
\end{eqnarray}
From the given distribution in Theorem~\ref{mainth}, we have
\begin{eqnarray}
&&\frac{P_{n_{J+1},...,n_{J+k}+1,...,n_{2J}}}
{P_{n_{(J+1,2J)}}}
=\frac{\lambda_{J+k}}{\kappa(n_{J+k}+1)\mu_{J+k}}
\elabel{dI}\\
&&\frac{P_{n_{J+1},...,n_{J+k}-1,...,n_{2J}}}
{P_{n_{(J+1,2J)}}}
=\frac{\kappa(n_{J+k})\mu_{J+k}}{\lambda_{J+k}}
\elabel{dII}\\
&&\frac{P_{n_{J+1},...,n_{J+r}+1,...,
n_{J+s}-1,...,n_{2J}}}
{P_{n_{(J+1,2J)}}}
=\frac{\kappa(n_{J+s})\lambda_{J+r}\mu_{J+s}}
{\kappa(n_{J+r}+1)\lambda_{J+s}\mu_{J+r}}.
\elabel{dIII}
\end{eqnarray}
Then, by \eq{balanck} to \eq{dIII}, we know
that it is sufficient to show
\begin{eqnarray}
&&\sum_{k=1}^{J}\alpha_{J+k}
+\sum_{k=1}^{J}\kappa(n_{J+k})\mu_{J+k}
\elabel{balancr}\\
&&=\sum_{k=1}^{J}
\left(1-\sum_{l=1}^{J}p_{J+k,J+l}\right)
\lambda_{J+k}+\sum_{k=1}^{J}\kappa(n_{J+k})
\alpha_{J+k}\mu_{J+k}/\lambda_{J+k}
\nonumber\\
&&\;\;\;\;+\sum_{r,s=1}^{J}
\kappa(n_{J+s})\lambda_{J+r}\mu_{J+s}
p_{{J+r},{J+s}}/\lambda_{J+s}.
\nonumber
\end{eqnarray}
Owing to the routing structure and the traffic equation \eq{traequ},
we have the following observations:
\begin{eqnarray}
&&\sum_{k=1}^{J}\left(1-\sum_{l=1}^{J}p_{J+k,J+l}
\right)\lambda_{J+k}=\sum_{k=1}^{J}\alpha_{J+k},
\elabel{hsum}\\
&&\sum_{r,s=1}^{J}
\kappa(n_{J+s})\lambda_{J+r}\mu_{J+s}
p_{{J+r},{J+s}}/\lambda_{J+s}
\elabel{dsum}\\
&&=\sum_{s=1}^{J}\left(\kappa(n_{J+s})\mu_{J+s}/
\lambda_{J+s}\right)
\sum_{r=1}^{J}\lambda_{J+r}p_{{J+r},{J+s}}
\nonumber\\
&&=\sum_{s=1}^{J}\kappa(n_{J+s})\mu_{J+s}
-\sum_{s=1}^{J}\kappa(n_{J+s})\mu_{J+s}
\alpha_{J+s}/\lambda_{J+s}.
\nonumber
\end{eqnarray}
Then, substituting \eq{hsum} and \eq{dsum} into
\eq{balancr}, we get the necessary equality.

Next, we derive the initial distribution $P^{s(J+k)}_{0}$
corresponding to state $n_{J+k}=0$. By the network stability
conditions \eq{kvp} and $\rho_{j}<1$ for each $j\in\{1,...,J\}$, it
follows from $\sum_{n_{J+k}=0}^{\infty}P^{s(J+k)}_{n_{J+k}}=1$ that
the following initial distribution is well posed,
\begin{eqnarray}
P^{s(J+k)}_{0}
=\frac{1}{1+\sum_{n_{J+k}=1}^{\infty}
\prod_{r_{J+k}=1}^{n_{J+k}}
\frac{\lambda_{J+k}}
{\kappa(r_{J+k})\mu_{J+k}}}.
\nonumber
\end{eqnarray}
Hence, we complete the proof of Theorem~\ref{mainth} for priority
types $h\in\{1,2\}$.

{\bf Step Three.} To finish the induction procedure,
in this step, we first suppose that the result
described in Theorem~\ref{mainth} is true for all
classes with priority rank $h\in\{1,...,I-1\}$, then,
we show that it is true for all
classes with priority rank $h\in\{1,...,I\}$.
By the similar illustration in getting \eq{i2de}, we
have the following differential equations,
\begin{eqnarray}
&&\dot{P}_{n_{((I-1)J+1,IJ)}}(t)
\elabel{ijde}\\
&&=-\sum_{k=1}^{J}\alpha_{(I-1)J+k}
P_{n_{((I-1)J,IJ)}}(t)
\nonumber\\
&&\;\;\;\;-\sum_{k=1}^{J}
\sum_{n_{1},...,n_{(I-1)J}=0}^{\infty}
a\left(n_{IJ+k},\sum_{u=0}^{I-2}n_{uJ+k}\right)
\mu_{(I-1)J+k}\left(1-p_{(I-1)J+k,(I-1)J+k}\right)
\nonumber\\
&&\;\;\;\;\;\;\;\;\;\;\;\;\;\;\;\;\;
\;\;\;\;\;\;\;\;\;\;\;\;\;\;\;\;\;
P_{n_{(1,(I-1)J)}}(t)P_{n_{((I-1)J+1,IJ)}}(t)
\nonumber\\
&&\;\;\;\;+\sum_{k=1}^{J}
\sum_{n_{1},...,n_{(I-1)J}=0}^{\infty}
a\left(n_{(I-1)J+k}+1,\sum_{u=0}^{I-2}
n_{uJ+k}\right)\mu_{(I-1)J+k}
\nonumber\\
&&\;\;\;\;\;\;\;\;\;\;\;\;\;\;
\left(1-\sum_{l=1}^{J}p_{(I-1)J+k,(I-1)J+l}\right)
P_{n_{(1,(I-1)J)}}(t)
P_{n_{(I-1)J+1},...,n_{(I-1)J+k}+1,...,n_{IJ}}(t)
\nonumber\\
&&\;\;\;\;+\sum_{k=1}^{J}
\left(n_{(I-1)J+k}\wedge 1\right)\alpha_{(I-1)J+k}
P_{n_{(I-1)J+1},...,n_{(I-1)J+k}-1,...,n_{IJ}}(t)
\nonumber\\
&&\;\;\;\;+\sum_{r,s=1,s\neq r}^{J}
\sum_{n_{1},...,n_{(I-1)J}=0}^{\infty}
a\left(n_{(I-1)J+r}+1,\sum_{u=0}^{I-2}n_{uJ+r}\right)
\mu_{(I-1)J+r}p_{(I-1)J+r,(I-1)J+s}
\nonumber\\
&&\;\;\;\;\;\;\;\;\;\;\;\;\;\;\;
P_{n_{(1,(I-1)J)}}(t)
\left(n_{(I-1)J+s}\wedge 1\right)
P_{n_{(I-1)J+1},...,n_{(I-1)J+r}+1,...,
n_{(I-1)J+s}-1,...,n_{IJ}}(t).
\nonumber
\end{eqnarray}
Next, we show that the given distribution corresponding the lowest
priority type in Theorem~\ref{mainth} is a steady-state solution of
those equations described by \eq{ijde}. It suffices to demonstrate
that the derivatives in the above equations are all made zero by
setting $P_{n_{((I-1)J,IJ)}}(t)=P_{n_{((I-1)J,IJ)}}$, i.e., to prove
that
\begin{eqnarray}
&&\sum_{k=1}^{J}\alpha_{(I-1)J+k}
P_{n_{((I-1)J,IJ)}}
\elabel{balde}\\
&&\;\;\;\;+\sum_{k=1}^{J}
\sum_{n_{1},...,n_{(I-1)J}=0}^{\infty}
a\left(n_{(I-1)J+k},\sum_{u=0}^{I-2}n_{uJ+k}\right)
P_{n_{(1,(I-1)J)}}\mu_{(I-1)J+k}P_{n_{((I-1)J+1,IJ)}}
\nonumber\\
&&=\sum_{k=1}^{J}
\sum_{n_{1},...,n_{(I-1)J}=0}^{\infty}
a\left(n_{(I-1)J+k}+1,\sum_{u=0}^{I-2}
n_{uJ+k}\right)P_{n_{(1,(I-1)J)}}
\nonumber\\
&&\;\;\;\;\;\;\;\;\;
\mu_{(I-1)J+k}
\left(1-\sum_{l=1}^{J}p_{(I-1)J+k,(I-1)J+l}\right)
P_{n_{(I-1)J+1},...,n_{(I-1)J+k}+1,...,n_{IJ}}
\nonumber\\
&&\;\;\;\;+\sum_{k=1}^{J}
\left(n_{(I-1)J+k}\wedge 1\right)\alpha_{(I-1)J+k}
P_{n_{(I-1)J+1},...,n_{(I-1)J+k}-1,...,n_{IJ}}
\nonumber\\
&&\;\;\;\;+\sum_{r,s=1,s\neq r}^{J}
\sum_{n_{1},...,n_{(I-1)J}=0}^{\infty}
a\left(n_{(I-1)J+r}+1,\sum_{u=0}^{I-2}n_{uJ+r}\right)
P_{n_{(1,(I-1)J)}}\mu_{(I-1)J+r}
\nonumber\\
&&\;\;\;\;\;\;\;\;\;
p_{(I-1)J+r,(I-1)J+s}\left(n_{(I-1)J+s}\wedge 1\right)
P_{n_{(I-1)J+1},...,n_{(I-1)J+r}+1,...,
n_{(I-1)J+s}-1,...,n_{IJ}}
\nonumber\\
&&\;\;\;\;+\sum_{k=1}^{J}
\sum_{n_{1},...,n_{(I-1)J}=0}^{\infty}
a\left(n_{(I-1)J+k},\sum_{u=0}^{I-2}n_{uJ+k}\right)
P_{n_{(1,(I-1)J)}}
\nonumber\\
&&\;\;\;\;\;\;\;\;\;\;\;\;\;\;\;\;\;\;\;
\;\;\;\;\;\;\;\;\;\;\;\;\;\;\;
\mu_{(I-1)J+k}p_{(I-1)J+k,(I-1)J+k}
P_{n_{((I-1)J+1,IJ)}}.
\nonumber
\end{eqnarray}
By the definition of $\kappa(\cdot)$ in Theorem~\ref{mainth}, we
have
\begin{eqnarray}
\kappa(n_{(I-1)J+k})=
\sum_{n_{1},...,n_{(I-1)J}=0}^{\infty}
a\left(n_{(I-1)J+k},\sum_{u=0}^{I-2}n_{uJ+k}\right)
P_{n_{(1,(I-1)J)}}.
\nonumber
\end{eqnarray}
Then, we can rewrite \eq{balde} as follows,
\begin{eqnarray}
&&\left(\sum_{k=1}^{J}\alpha_{(I-1)J+k}
+\sum_{k=1}^{J}\kappa(n_{(I-1)J+k})
\mu_{(I-1)J+k}\right)P_{n_{((I-1)J+1,IJ)}}
\elabel{balancK}\\
&&=\sum_{k=1}^{J}\kappa(n_{(I-1)J+k}+1)\mu_{(I-1)J+k}
\left(1-\sum_{l=1}^{J}p_{(I-1)J+k,(I-1)J+l}\right)
\nonumber\\
&&\;\;\;\;\;\;\;\;\;
P_{n_{(I-1)J+1},...,n_{(I-1)J+k}+1,...,n_{IJ}}
\nonumber\\
&&\;\;\;\;+\sum_{k=1}^{J}
\left(n_{(I-1)J+k}\wedge 1\right)\alpha_{(I-1)J+k}
P_{n_{(I-1)J+1},...,n_{(I-1)J+k}-1,...,n_{IJ}}
\nonumber\\
&&\;\;\;\;+\sum_{r,s=1,s\neq r}^{J}
\kappa(n_{(I-1)J+k}+1)
\mu_{(I-1)J+r}p_{(I-1)J+r,(I-1)J+s}
\left(n_{(I-1)J+s}\wedge 1\right)
\nonumber\\
&&\;\;\;\;\;\;\;\;\;\;\;\;\;\;\;\;\;\;\;\;
P_{n_{(I-1)J+1},...,n_{(I-1)J+r}+1,...,
n_{(I-1)J+s}-1,...,n_{IJ}}
\nonumber\\
&&\;\;\;\;+\sum_{k=1}^{J}\kappa(n_{(I-1)J+k})
\mu_{(I-1)J+k}p_{(I-1)J+k,(I-1)J+k}
P_{n_{((I-1)J+1,IJ)}}.
\nonumber
\end{eqnarray}
From the given distribution in Theorem~\ref{mainth}, we have
\begin{eqnarray}
&&\frac{P_{n_{(I-1)J+1},...,n_{(I-1)J+k}+1,...,
n_{2J}}}{P_{n_{((I-1)J+1,2J)}}}
=\frac{\lambda_{(I-1)J+k}}
{\kappa(n_{(I-1)J+k}+1)\mu_{(I-1)J+k}}
\elabel{dIK}\\
&&\frac{P_{n_{(I-1)J+1},...,n_{(I-1)J+k}-1,...,
n_{2J}}}{P_{n_{((I-1)J+1,2J)}}}
=\frac{\kappa(n_{(I-1)J+k})
\mu_{(I-1)J+k}}{\lambda_{(I-1)J+k}}
\elabel{dIIK}\\
&&\frac{P_{n_{(I-1)J+1},...,n_{(I-1)J+r}+1,...,
n_{(I-1)J+s}-1,...,n_{2J}}}
{P_{n_{((I-1)J+1,2J)}}}
=\frac{\kappa(n_{(I-1)J+s})\lambda_{(I-1)J+r}
\mu_{(I-1)J+s}}{\kappa(n_{(I-1)J+r}+1)
\lambda_{(I-1)J+s}\mu_{(I-1)J+r}}.
\elabel{dIIIK}
\end{eqnarray}
From \eq{balancK} to \eq{dIIIK}, it is sufficient to show
\begin{eqnarray}
&&\sum_{k=1}^{J}\alpha_{(I-1)J+k}
+\sum_{k=1}^{J}\kappa(n_{(I-1)J+k})\mu_{(I-1)J+k}
\elabel{balancrK}\\
&&=\sum_{k=1}^{J}
\left(1-\sum_{l=1}^{J}p_{(I-1)J+k,(I-1)J+l}\right)
\lambda_{(I-1)J+k}
\nonumber\\
&&\;\;\;\;+\sum_{k=1}^{J}\kappa(n_{(I-1)J+k})
\alpha_{(I-1)J+k}\mu_{(I-1)J+k}/\lambda_{(I-1)J+k}
\nonumber\\
&&\;\;\;\;+\sum_{r,s=1}^{J}
\kappa(n_{(I-1)J+s})\lambda_{(I-1)J+r}\mu_{(I-1)J+s}
p_{{(I-1)J+r},{(I-1)J+s}}/\lambda_{(I-1)J+s}.
\nonumber
\end{eqnarray}
Owing to the routing structure and the traffic equation \eq{traequ},
we have the following observations:
\begin{eqnarray}
&&\sum_{k=1}^{J}\left(1-\sum_{l=1}^{J}
p_{(I-1)J+k,(I-1)J+l}\right)\lambda_{(I-1)J+k}
=\sum_{k=1}^{J}\alpha_{(I-1)J+k},
\elabel{hsumK}\\
&&\sum_{r,s=1}^{J}
\kappa(n_{(I-1)J+s})\lambda_{(I-1)J+r}\mu_{(I-1)J+s}
p_{{(I-1)J+r},{(I-1)J+s}}/\lambda_{(I-1)J+s}
\elabel{dsumK}\\
&&=\sum_{s=1}^{J}
\left(\kappa(n_{(I-1)J+s})\mu_{(I-1)J+s}/
\lambda_{(I-1)J+s}\right)
\sum_{r=1}^{J}\lambda_{(I-1)J+r}
p_{{(I-1)J+r},{(I-1)J+s}}
\nonumber\\
&&=\sum_{s=1}^{J}\kappa(n_{(I-1)J+s})
\mu_{(I-1)J+s}
-\sum_{s=1}^{J}\kappa(n_{(I-1)J+s})\mu_{(I-1)J+s}
\alpha_{(I-1)J+s}/\lambda_{(I-1)J+s}.
\nonumber
\end{eqnarray}
Then, substituting \eq{hsumK} and \eq{dsumK} into
\eq{balancrK}, we get the necessary equality.

Next, we derive the initial distribution $P^{s((I-1)J+k)}_{0}$
corresponding to state $n_{(I-1)J+k}=0$. By the network stability
conditions \eq{kvp} and $\rho_{j}<1$ for each $j\in\{1,...,J\}$, it
follows from $\sum_{n_{(I-1)J+k}=0}^{\infty}
P^{s((I-1)J+k)}_{n_{(I-1)J+k}}=1$ that the following initial
distribution is well posed,
\begin{eqnarray}
P^{s((I-1)J+k)}_{0}
=\frac{1}{1+\sum_{n_{(I-1)J+k}=1}^{\infty}
\prod_{r_{(I-1)J+k}=1}^{n_{(I-1)J+k}}
\frac{\lambda_{(I-1)J+k}}
{\kappa(r_{(I-1)J+k})\mu_{(I-1)J+k}}}.
\nonumber
\end{eqnarray}
Hence, we complete the proof of Theorem \ref{mainth}. $\Box$

\subsection{Proof of Proposition~\ref{sn}}

\indent {\bf Net I.} We justify the stability condition \eq{kvp} for
Net I by induction in terms of the number of job types, i.e.,
$h=1,2,...,I$.

First, we consider the case that $h=2$. Note that $c_{j}=1$ for all
$j\in\{1,2,...,J\}$ and step one in the proof of
Theorem~\ref{mainth}, we have that, for $n_{J+k}>0$ with
$k\in\{1,2,...,J\}$,
\begin{eqnarray}
\kappa(n_{J+k})=P_{0}^{s(k)}=1-\lambda_{k}/\mu_{k}.
\nonumber
\end{eqnarray}
Therefore, the condition \eq{kvp} is true for $h=2$. Furthermore, we
have
\begin{eqnarray}
P^{s(J+k)}_{0}=1-\frac{\lambda_{J+k}}
{\left(1-\lambda_{k}/\mu_{k}\right)\mu_{J+k}}.
\nonumber
\end{eqnarray}

Second, for $h\leq I-1$, we suppose that
\begin{eqnarray}
\kappa(n_{(h-1)J+k}) =1-\sum_{i=0}^{h-2}\lambda_{iJ+k}/\mu_{iJ+k}.
\elabel{c11}
\end{eqnarray}
Hence, we have
\begin{eqnarray}
P^{s((h-1)J+k)}_{0}=1-\frac{\lambda_{(h-1)J+k}}
{\left(1-\sum_{i=0}^{h-2}\lambda_{iJ+k}/\mu_{iJ+k}
\right)\mu_{(h-1)J+k}}.
\elabel{c12}
\end{eqnarray}
From the induction assumptions \eq{c11} and \eq{c12}, we know that
\eq{kvp} is true for $h\leq I-1$.

Finally, we show that \eq{kvp} holds for
$h=I$. In fact, from the definition of
$\kappa(n_{(I-1)J+k})$, we know that, for
$n_{(I-1)J+k}>0$,
\begin{eqnarray}
\kappa(n_{(I-1)J+k})=\prod_{i=0}^{I-2}P^{s(iJ+k)}_{0}
=1-\sum_{i=0}^{I-2}\lambda_{iJ+k}/\mu_{iJ+k}.
\nonumber
\end{eqnarray}
Then, we see that \eq{kvp} is true for $h=I$. Hence, we complete the
proof of Net I.

\noindent {\bf Net II.}  Note that, for each $k\in\{1,...,J\}$ and
$n_{J+k}\geq c_{s(k)}$, we have
\begin{eqnarray}
&&\kappa(n_{J+k})
\nonumber\\
&&=\sum_{n_{1},...,n_{J}=0}^{\infty}
a(n_{J+k},n_{k})P_{n_{(1,J)}}
\nonumber\\
&&=\sum_{n_{k}=0}^{c_{s(k)}-1}
a(n_{J+k},n_{k})P^{s(k)}_{n_{k}}
\nonumber\\
&&=P^{s(k)}_{0}\left(c_{s(k)}
+(c_{s(k)}-1)(\lambda_{k}/\mu_{k})
+...+
\frac{(\lambda_{k}/\mu_{k})^{c_{s(k)}-1}}
{(c_{s(k)}-1)!}\right)
\nonumber\\
&&=\frac{\sum_{n_{k}=0}^{c_{s(k)}-1}
(c_{s(k)}-n_{k})
\frac{(\lambda_{k}/\mu_{k})^{n_{k}}}{n_{k}!}}
{\sum_{n_{k}=0}^{c_{s(k)}-1}
\frac{(\lambda_{k}/\mu_{k})^{n_{k}}}{n_{k}!}
+\frac{\left(\lambda_{k}/\mu_{k}\right)^{c_{s(k)}}}
{c_{s(k)}!\left(1-\lambda_{k}/c_{s(k)}\mu_{k}\right)}}
\nonumber\\
&&=
\frac{c_{s(k)}\left(1-\lambda_{k}/c_{s(k)}\mu_{k}\right)
\sum_{n_{k}=0}^{c_{s(k)}-1}
\frac{\left(c_{s(k)}-n_{k}\right)
\left(\lambda_{k}/\mu_{k}\right)^{n_{k}}}
{c_{s(k)}n_{k}!}}
{1+\sum_{n_{k}=0}^{c_{s(k)}-2}
\left(\frac{1}{(n_{k}+1)!}-
\frac{1}{c_{s(k)}n_{k}!}\right)
(\frac{\lambda_{k}}{\mu_{k}})^{n_{k}+1}}
\nonumber\\
&&= c_{s(k)}\left(1-\lambda_{k}/c_{s(k)}\mu_{k}\right). \nonumber
\end{eqnarray}
Furthermore, for $0<n_{J+k}<c_{s(k)}$, we have
\begin{eqnarray}
&&\kappa(n_{J+k})
\nonumber\\
&&=\sum_{n_{1},...,n_{J}=0}^{\infty}
a(n_{J+k},n_{k})P_{n_{(1,J)}}
\nonumber\\
&&=\sum_{n_{k}=0}^{c_{s(k)}-1}
a(n_{J+k},n_{k})P^{s(k)}_{n_{k}}
\nonumber\\
&&=P_{0}^{s(k)}\left(n_{J+k}
\left(1+\frac{1}{1!}
\left(\frac{\lambda_{k}}{\mu_{k}}\right)
+...+\frac{1}{(c_{s(k)}-n_{J+k})!}
\left(\frac{\lambda_{k}}{\mu_{k}}
\right)^{c_{s(k)}-n_{J+k}}\right)\right.
\nonumber\\
&&\;\;\;\;\;\;\;\;\;\;\;\;\;\;\;
+\left.\frac{n_{J+k}-1}{(c_{s(k)}-n_{J+k}+1)!}
\left(\frac{\lambda_{k}}{\mu_{k}}
\right)^{c_{s(k)}-n_{J+k}+1}+...+
\frac{1}{(c_{s(k)}-1)!}
\left(\frac{\lambda_{k}}{\mu_{k}}
\right)^{c_{s(k)}-1}\right)
\nonumber\\
&&=
\frac{n_{J+k}\sum_{n_{k}=0}^{c_{s(k)}-n_{J+k}}
\frac{\left(\lambda_{k}/\mu_{k}
\right)^{n_{k}}}{n_{k}!}+
\sum_{n_{k}=c_{s(k)}-n_{J+k}+1}^{c_{s(k)}-1}
\frac{\left(c_{s(k)}-n_{k}\right)
\left(\lambda_{k}/\mu_{k}\right)^{n_{k}}}{n_{k}!}}
{\sum_{n_{k}=0}^{c_{s(k)}-1}
\frac{(\lambda_{k}/\mu_{k})^{n_{k}}}{n_{k}!}
+\frac{\left(\lambda_{k}/\mu_{k}\right)^{c_{s(k)}}}
{c_{s(k)}!\left(1-\lambda_{k}/c_{s(k)}\mu_{k}\right)}}
\nonumber\\
&&=
\frac{n_{J+k}\left(1-\lambda_{k}/c_{s(k)}\mu_{k}\right)
\left(\sum_{n_{k}=0}^{c_{s(k)}-n_{J+k}}
\frac{\left(\lambda_{k}/\mu_{k}
\right)^{n_{k}}}{n_{k}!}+
\sum_{n_{k}=c_{s(k)}-n_{J+k}+1}^{c_{s(k)}-1}
\frac{\left(c_{s(k)}-n_{k}\right)
\left(\lambda_{k}/\mu_{k}\right)^{n_{k}}}
{n_{k}!n_{J+k}}\right)}
{1+\sum_{n_{k}=0}^{c_{s(k)}-2}
\left(\frac{1}{(n_{k}+1)!}-
\frac{1}{c_{s(k)}n_{k}!}\right)
(\lambda_{k}/\mu_{k})^{n_{k}+1}}
\nonumber\\
&&\geq n_{J+k}
\left(1-\lambda_{k}/c_{s(k)}\mu_{k}\right).
\nonumber
\end{eqnarray}
Then, we can see that \eq{kvp} is true. Hence, we complete the proof
of Proposition~\ref{sn}. $\Box$

\subsection{Illustration of
Example~\ref{exampleI}}

For each class $k\in\{1,2,...,J\}$, the initial distribution
$P_{0}^{s(k)}$ can be calculated as follows since $\rho_{s(k)}<1$,
\begin{eqnarray}
P_{0}^{s(k)}=\frac{1-\lambda_{k}/2\mu_{k}} {1+\lambda_{k}/2\mu_{k}}.
\nonumber
\end{eqnarray}
Then, for $n_{J+k}=1$, we have
\begin{eqnarray}
\kappa(n_{J+k})=P_{0}^{s(k)} \left(1+\lambda_{k}/\mu_{k}\right)
=\frac{(1-\lambda_{k}/2\mu_{k})(1+\lambda_{k}/\mu_{k})}
{1+\lambda_{k}/2\mu_{k}}\geq 1-\lambda_{k}/2\mu_{k}. \nonumber
\end{eqnarray}
Furthermore, for $n_{J+k}\geq 2$, we have
\begin{eqnarray}
\kappa(n_{J+k})=P_{0}^{s(k)}(2+\lambda_{k}/\mu_{k})
=2(1-\lambda_{k}/2\mu_{k}). \nonumber
\end{eqnarray}
Therefore, since $\rho_{s(k)}<1$, we have
\begin{eqnarray}
P^{s(J+k)}_{0} &=&\frac{1}{1+\frac{\lambda_{J+k}}
{\kappa(1)\mu_{J+k}} +\frac{\lambda_{J+k}}{\kappa(1)\mu_{J+k}}
\frac{\lambda_{J+k}/\mu_{J+k}} {(1-\lambda_{J+k}/\kappa(2)\mu_{J+k})
\kappa(2)\mu_{J+k}}}
\nonumber\\
&=& \frac{\left(1-\lambda_{k}/2\mu_{k}-
\lambda_{J+k}/2\mu_{J+k}\right) \left(1+\lambda_{k}/\mu_{k}\right)}
{1+\lambda_{k}/2\mu_{k}+\lambda_{J+k}/2\mu_{J+k}
-\frac{1}{2}(\lambda_{k}/\mu_{k})^{2}}. \nonumber
\end{eqnarray}
Finally, we can compute $\kappa(n_{2J+k})$ for $n_{2J+k}\geq 2$ as
follows,
\begin{eqnarray}
&&\kappa(n_{2J+k})
\nonumber\\
&&= 2P_{0}^{s(k)}P_{0}^{s(J+k)}+ P_{1}^{s(k)}P_{0}^{s(J+k)}
+P_{0}^{s(k)}P_{1}^{s(J+k)}
\nonumber\\
&&= P_{0}^{s(k)}P_{0}^{s(J+k)} \left(2+\lambda_{k}/\mu_{k}+
\lambda_{J+k}/\mu_{J+k}\right)
\nonumber\\
&&= \frac{2\left(1-\frac{\lambda_{k}}{2\mu_{k}}-
\frac{\lambda_{J+k}}{2\mu_{J+k}}\right)
\left(1+\frac{\lambda_{k}}{\mu_{k}}\right)
\left(1-\frac{\lambda_{k}}{2\mu_{k}}\right)
\left(1+\frac{\lambda_{k}}{2\mu_{k}}+
\frac{\lambda_{J+k}}{2\mu_{J+k}}\right)}
{\left(1+\frac{\lambda_{k}}{2\mu_{k}}+
\frac{\lambda_{J+k}}{2\mu_{J+k}}
-\frac{1}{2}\left(\frac{\lambda_{k}} {\mu_{k}}\right)^{2}\right)
\left(1+\frac{\lambda_{k}}{2\mu_{k}}\right)}
\nonumber\\
&&= \frac{2\left(1-\frac{\lambda_{k}}{2\mu_{k}}-
\frac{\lambda_{J+k}}{2\mu_{J+k}}\right)
\left(1+\frac{\lambda_{k}}{\mu_{k}}
+\frac{\lambda_{J+k}}{2\mu_{J+k}}
-\left(\frac{\lambda_{k}}{2\mu_{k}}\right)^{2} -\frac{1}{4}
\left(\frac{\lambda_{k}}{\mu_{k}}\right)^{3}
+2\left(1-\frac{\lambda_{k}}{2\mu_{k}}\right)
\frac{\lambda_{k}\lambda_{J+k}}{4\mu_{k}\mu_{J+k}}\right)}
{1+\frac{\lambda_{k}}{\mu_{k}} +\frac{\lambda_{J+k}}{2\mu_{J+k}}
-\left(\frac{\lambda_{k}}{2\mu_{k}}\right)^{2} -\frac{1}{4}
\left(\frac{\lambda_{k}}{\mu_{k}}\right)^{3}
+\frac{\lambda_{k}\lambda_{J+k}}{4\mu_{k}\mu_{J+k}}}. \nonumber
\end{eqnarray}
Thus, if $\lambda_{k}/\mu_{k}\leq 1$, from the condition
$\rho_{s(k)}<1$, we see that
\begin{eqnarray}
\frac{\lambda_{2J+k}}{\kappa(n_{2J+k})\mu_{2J+k}}
\leq\frac{\lambda_{2J+k}} {2\left(1-\lambda_{k}/2\mu_{k}
-\lambda_{J+k}/2\mu_{J+k}\right)\mu_{2J+k}}<1, \nonumber
\end{eqnarray}
that is, the condition \eq{kvp} holds. $\Box$

\section{Conclusion}\label{conclusion}

The research conducted in the paper is to iteratively derive the
product-form solutions of stationary distributions of priority
multiclass queueing networks with multi-sever stations. The networks
are Markovian with exponential interarrival and service time
distributions. These solutions can be used to conduct performance
analysis or as comparison criteria for approximation and simulation
studies of large scale networks with multi-processor shared-memory
switches and cloud computing systems with parallel-server stations.
Numerical comparisons with existing Brownian approximating model are
provided to show the effectiveness of our algorithm. For more
discussions related to general distributions of interarrival and
service times, performance analysis and optimization, readers are
referred to Dai {\em et
al.}~\cite{dai:broapp,dai:difapp,daidai:heatra,dai:optrat}.

\bibliographystyle{abbrv}
\bibliography{dai}
\end{document}